\documentclass{m2an}
\usepackage{amsmath,amssymb,amsfonts,setspace,bbold}
\usepackage{lmodern}
\usepackage{graphicx,caption,epsfig,subfigure,epsfig,wrapfig,sidecap}
\usepackage{url,color,xcolor,verbatim,algorithmicx,dsfont}

\def\R{\mathbb{R}}

\def\E{\mathbb{E}}

\def\rd{\mathrm{d}}

\definecolor{mygrey}{gray}{0.75}

\newenvironment{rmat}{\left[\begin{array}{rrrrrrrrrrrrr}}{\end{array}\right]}
\newcommand\brm{\begin{rmat}}
\newcommand\erm{\end{rmat}}
\newenvironment{cmat}{\left[\begin{array}{ccccccccc}}{\end{array}\right]}
\newcommand\bcm{\begin{cmat}}
\newcommand\ecm{\end{cmat}}

\numberwithin{equation}{section}
\begin{document}
%%-----------------------------
%%      the top matter
%%-----------------------------
\title{Numerical stability revisited: A family of benchmark problems for the analysis of explicit stochastic differential equation integrators
}\thanks{T. Hudson is grateful for support from UKRI during this work through NSF-UKRI grant reference UKRI3611.}\thanks{ X.~Li and S. Helfert are grateful for partial support by the NSF Award DMS-1847770, NSF-UKRI DMS-2513795 and the 2023 UNC Charlotte faculty research grant.}\thanks{The authors would like to thank the International Centre for the Mathematical Sciences, Edinburgh, for hosting them for a week-long workshop during the course of this work.}% At most 5 thanks

% Authors: full names plus addresses.
\author{Thomas Hudson}\address{Warwick Mathematics Institute, University of Warwick, Coventry, UK; Email: {t.hudson.1@warwick.ac.uk}, (\url{tinyurl.com/thudso}).}
\author{Sarah Helfert}\address{Department of  Mathematics \& Statistics, University of North Carolina Charlotte, NC, USA; Email: {shelfert@charlotte.edu}.}
\author{Xingjie Helen Li}\address{Department of Mathematics \& Statistics, University of North Carolina Charlotte, NC, USA; Email: {xli47@charlotte.edu}.}

%% abstract 
\begin{abstract}
We revisit the numerical stability of four well-established explicit stochastic integration schemes through a new generic benchmark stochastic differential equation designed to assess asymptotic statistical accuracy and stability properties. This one-parameter benchmark equation is derived from a general one-dimensional first-order SDE using spatio-temporal nondimensionalization and is employed to evaluate the performance of the (1) Euler–Maruyama, (2) Milstein, (3) Stochastic Heun, and (4) three-stage Runge–Kutta schemes.
Our findings reveal that lower-order schemes can outperform higher-order ones over a range of time step sizes, depending on the benchmark parameters and application context. The theoretical results are validated through a series of numerical experiments, and we discuss their implications for more general applications, including a nonlinear example. Our results suggest that the insights obtained from the linear benchmark problem provide reliable guidance for time-stepping strategies when simulating nonlinear SDEs.
\end{abstract}

\maketitle%%-----------------------------
%%      the top matter
%%-----------------------------
\title{Numerical stability revisited: A family of benchmark problems for the analysis of explicit stochastic differential equation integrators
}
% Authors: full names plus addresses.
\author{Thomas Hudson \thanks{Mathematics Institute, University of Warwick, Coventry, CV4 7AL, UK 
  ({t.hudson.1@warwick.ac.uk}). TH is grateful for support from from UKRI during this work through NSF-UKRI grant reference UKRI3611.}
  \and Sarah Helfert \thanks{Department of  Mathematics \& Statistics, University of North Carolina Charlotte, NC 22833,   ({shelfert@charlotte.edu}).  S. Murphy is grateful for partial support by the NSF Award DMS-1847770 and NSF-UKRI DMS-2513795.}
\and Xingjie Helen Li \thanks{Department of Mathematics \& Statistics, University of North Carolina Charlotte, NC 22833, 
  ({xli47@charlotte.edu}). X.~Li is grateful for partial support by the NSF Award DMS-1847770, NSF-UKRI DMS-2513795 and the 2023 UNC Charlotte faculty research grant.}
}

%\subjclass{60H35, 65L20}
%\keywords{Analysis of explicit numerical integrators for SDEs; asymptotic statistical stability; spatio-temporal nondimensionalization}
\maketitle

\section{Introduction}
Stochastic differential equation (SDE) models are important for modelling a wide range of real--world phenomena, and combine both deterministic and random effects into a model describing the evolution of degrees of freedom in time. In many cases, SDE models are ergodic, entailing that the law of the process converges to an equilibrium distribution independently of the initial condition. In practice, ergodicity is a desirable property for real-world applications, with the equilibrium distribution representing a form of thermalised steady-state for a system.

To simulate and predict the behaviour of SDE models, numerical integration algorithms are needed, and a wide range of algorithms have therefore been developed to simulate SDE models. Even when an SDE is ergodic, it has long been recognized that a numerical discretisation of the evolution may not reproduce the equilibrium distribution, leading to biased long-time statistics. For example, when the drift and diffusion coefficients are only \textit{locally Lipschitz}, explicit schemes can fail to be ergodic and therefore yield biased asymptotic statistics, even when the underlying SDE is itself geometrically ergodic \cite{robert1996exponential,mattingly2002ergodicity,fang2020adaptive,kelly2022adaptive}. On the other hand, when the coefficients are globally Lipschitz, explicit schemes admit strong and weak convergence guarantees and can be proved to reproduce the equilibrium distribution when time step size $h$ tends to $0$  \cite{saito1991discrete,kloeden2012numerical,lord2014stochastic}. Moreover, for a fixed simulation time $T$, higher-order schemes are provably more accurate than lower-order schemes as the time step tends to zero. However, the associated error bounds typically involve prefactors that depend on $T$, so these results do not by themselves imply that a higher-order method will consistently outperform a lower-order one in the long-time limit for a fixed time-step size. Motivated by this, we are interested here in comparing the long-time behavior of explicit schemes of a range of orders in the regime where $T\to\infty$ while the time step $h$ remains moderately large, thereby focusing on the ability of the schemes we consider to approximate the invariant measure and the resulting asymptotic statistics. 

We note at this stage that where discretisations of SDEs are used as a method to generate samples of a particular distribution with a known density, Monte Carlo algorithms often use an accept-reject or Metropolis step to counteract the bias inherent in the numerical discretisation \cite{cances2007theoretical,legoll2022adaptive}. On the other hand, when we are interested not only in sampling the equilibrium distribution but also in the dynamics of the process, performing accept-reject steps on entire trajectory segments may be costly, as they involve discarding expensive-to-compute trajectories of high-dimensional systems of SDEs. As an alternative, we ask how accurately and stably numerical schemes are able to capture asymptotic statistical properties of the underlying SDE model. If we are able to quantify the errors committed, and these are at an acceptable level for the required application, we may as a result be able to save extra computational work. Furthermore, stochasticity and stiffness appear in a wide range of physical and chemical systems, and the time-step size is often severely restricted to maintain numerical stability and accuracy when integrating stiff stochastic
differential equations with fast diffusion process \cite{han2017numerical,han2019explicit,han2022modeling}.

Classically, the stability of SDE schemes has been measured through an approach known as mean-square stability analysis \cite{schurz1999preservation,saito1996stability,saito2002mean,tocino2012mean}. These approaches have tended to use geometric Brownian motion (GBM) as a benchmark, i.e.
\begin{equation*}
    \rd x_t = -x_t \rd t + \eta x_t\rd W_t,
\end{equation*}
extending the concept of linear stability from deterministic ODEs \cite{saito1996stability,kloeden2012numerical}.
In this work, we revisit the underlying assumptions of this stability analysis, and propose a more generic benchmark SDE problem for the study of statistical accuracy and stability properties, namely the one-dimensional equation
\begin{equation*}
\rd x_t = -x_t \rd t + (1+\eta x_t)\rd W_t.
\end{equation*}
Here, $\eta$ is a real parameter which controls the relative strength of the multiplicative noise and the stiffness of the deterministic drift. By non-dimensionalising, we show that this benchmark arises generically from a general first-order SDE in one dimension with affine coefficients.
While the study of the GBM problem benefits from a simple exact expression for solution trajectories in terms of Brownian motion, we will argue below that since it represents a less generic situation, it does not necessarily reflect the performance of numerical integration schemes in a realistic nonlinear setting. Despite their superficial similarity, as a stark illustration of the contrast between the two problems, the equilibrium distribution for GBM is always concentrated at $x=0$ when it exists, while our proposed benchmark includes a wide range of algebraic tail behaviour for the resulting equilibrium distribution.

At this stage, we briefly comment on our reasons for selecting a one-dimensional rather than a multi-dimensional test case. In higher dimensions, a non-dimensionalisation strategy similar to that which we pursue below to derive our benchmark will lead to a much richer parametric representation, which is therefore more complex to analyse. These challenges are already recognised in the literature on mean-squared stability; see for example \cite{schurz1996asymptotical,saito2002mean}, and further discussion the complexities faced in higher-dimensional cases is provided when deriving our benchmark below.

After deriving the benchmark problem, we use it to analyse the long-time performance of some well--established explicit integration schemes. Our focus on explicit schemes is motivated by the need for integrators which avoid costly implicit steps in high-dimensional problems. While there is no simple exact formula for solutions of our test problem, we will see that it is nevertheless simple enough to allow for the explicit computation of both the equilibrium distribution and the evolution of the low-order moments of the law as the process evolves. We choose the latter as metrics for assessing the accuracy of the schemes, reflecting the sorts of statistics which are often of interest to practitioners.

Along with a discussion of the {long-time} statistical accuracy of the numerical methods considered, we obtain results characterising the range of time-steps for which the methods are stable in the sense that they preserve the boundedness of first and second moments as the number of time steps simulated tends to infinity. This is a feature which we term \emph{statistical stability} to distinguish it from other differing notions \cite{saito1996stability,milstein2013numerical,kloeden2012numerical}. The ranges of time step we find provide us with a meaningful way to compare the stability of methods as the parameter $\eta$ varies, allowing us to represent a range of relative strengths for the multiplicative noise and drift. In modeling contexts where long trajectories are important in order to enable the observation of phenomena which would otherwise be out of reach, understanding {and} predicting the stability of methods are important as a way to allow for increased time-step sizes, which trade accuracy for physical insight which would otherwise not be available.

Finally, to verify the analytical results we obtain, we conclude by performing a range of numerical experiments, and discuss the consequences of our results for a more general application to a nonlinear example.\footnote{Accompanying codes can be found at \url{https://github.com/XingjieHelenLi/NumSDE_Revisit}.}

\medskip
\paragraph{\bf Outline.}
We now provide a brief guide to the remainder of the paper. In Section~\ref{sec:schemes}, we review the numerical schemes we study, and give an overview of established results on their local accuracy properties. In Section~\ref{sec:benchmark}, we use a nondimensionalisation argument to derive our test problem from a general linear SDE, discussing various exactly computable properties of the resulting equation. In Section~\ref{sec:stab}, we derive expressions for the evolution of discrete statistical moments under various numerical integrations, comparing to the continuous counterparts for the benchmark problem, verifying the results numerically and studying the corresponding stability regions. In Section~\ref{sec:nonlinear}, we demonstrate the utility of our results in the context of a more realistic setting of an autonomous equation where the drift and diffusion coefficients are both nonlinear.

\section{Overview of numerical schemes}
\label{sec:schemes}
We consider the stochastic initial value problem for the scalar
autonomous It\^{o} stochastic differential equation (SDE) given by
\begin{equation}\label{eq:Ito_SDE}
\rd x_t =f(x_t)\rd t + g(x_t)\rd W_t,\quad t\in[0,T].
\end{equation}
We will assume throughout that the solution is subject to a deterministic initial condition, $x_0$.

For the purposes of comparing the numerical schemes we study, we will discretise the simulation time interval into subintervals of fixed length $h = T/N$, and define a discrete time grid $t_n = nh$ for $n=0,\ldots,N$. For the different schemes we discuss, $x_n$ will always be used to refer to the approximation of the solution generated at time $t_n$. For convenience in describing the numerical schemes below, we also introduce the notation
\[
f_n := f(x_n),\qquad\text{and}\qquad f'_n:=f'(x_n),
\]
and use analogous expressions $g_n$, $g'_n$ and $g''_n$ for the relevant functions evaluated at $x_n$.
We also introduce sequences of independent identically distributed Gaussian random variables $\Delta W_n$ and $\Delta \tilde{W}_n$ with mean zero and variance $h$, i.e. $\Delta \tilde{W}_n,\Delta W_n\sim\mathcal{N}(0,h)$.

We select four classical explicit numerical schemes to study with a range of levels of accuracy. In particular, we consider:
\begin{enumerate}
    \item The Euler--Maruyama (EM) method \cite{Maruyama55}, \cite[\S9.1]{kloeden2012numerical}:
    \begin{equation}\label{eq:EM}
        x_{n+1}=x_n + f_n h+g_n \Delta W_n.     \tag{EM}
    \end{equation}
    \item The Milstein (Mil) method \cite[\S10.3]{kloeden2012numerical}:
    \begin{equation}\label{eq:Mil}
    x_{n+1}=x_n + f_n h+g_n \Delta W_n+\tfrac{1}{2}g_n'g_n \left(\Delta W_n^2-h\right).\tag{Mil}
    \end{equation}
     \item The Stochastic Heun (SH) method \cite{saito1996stability}:
     \begin{equation}\label{eq:Heu}
         x_{n+1} =x_{n}+\tfrac{1}{2}\left[F_1+F_2\right]h
+\tfrac{1}{2}\left[G_1+G_2\right]\Delta W_n,\tag{SH}
     \end{equation}
     where setting $F(x):= f(x)-\tfrac{1}{2}g'(x)g(x)$, the coefficients in the relation above are defined to be
\[
\begin{aligned}
F_1 &= F(x_n), &
G_1 &= g(x_n),\\
F_2 &= F(x_n+F_1h+G_1\Delta W_n),&
G_2 &= g(x_n+F_1h+G_1\Delta W_n).
\end{aligned}
\]
    \item The improved 3-stage Runge-Kutta (RK3) method \cite{saito1991discrete,saito1996stability}:
    \begin{align}\label{eq:RK3}
    %\begin{aligned}
    x_{n+1}&= x_n+\tfrac{1}{4}\left(F_1+3F_3\right)h
    +\tfrac{1}{4}\left[G_1+3G_3\right]\Delta W_n \nonumber \\
    &\qquad +\tfrac{1}{2\sqrt{3}}\left(f'_ng_n-g'_nf_n-\tfrac{1}{2}g''_ng^2_n\right)h\Delta \tilde{W}_n,\tag{RK3}
    %\end{aligned}
    \end{align}
    where we again denote $F(x):= f(x)-\tfrac{1}{2}g'(x)g(x)$, and the coefficients are
    \[
    \begin{aligned}
        F_1 &= F(x_n), &
        G_1 &= g(x_n),\\
        F_2 &= F(x_n+\tfrac{1}{3}F_1h+\tfrac{1}{3} G_1 \Delta W_n),& 
        G_2 &= g(x_n+\tfrac{1}{3}F_1h+\tfrac{1}{3} G_1 \Delta W_n),\\
        F_3 &= F(x_n+\tfrac{2}{3}F_2h+\tfrac{2}{3}G_2 \Delta W_n),& 
        G_3 &= g(x_n+\tfrac{2}{3}F_2h+\tfrac{2}{3}G_2 \Delta W_n).
    \end{aligned}
    \]
\end{enumerate}
Established results on the strong and weak orders of these schemes are summarised in Table~\ref{tab:str_weak_accuracy_summary}, and the accuracy in terms of time-stepping $h$ is summarised in Table~\ref{tab:accuracy_summary}. As we can see, these schemes range significantly across the convergence order that they exhibit.

\begin{table}[tp!]
    \centering
    \begin{tabular}{r|c|c}
        \textbf{Method} & \textbf{Strong order } & \textbf{Weak order} \\
        \hline
        Euler-Maruyama (EM) &  $0.5$ & $1$ \\
        Milstein (Mil) & $1$ & $1$ \\
        Stochastic Heun (SH) & $1$ & $2$\\
        3-stage Runge-Kutta (RK3) & $1$ & $2$
    \end{tabular}
    \caption{A summary of the strong and weak convergence orders of the methods studied. The first two rows correspond to results proved in \cite{kloeden2012numerical}, while the latter two rows correspond to results proved in \cite{saito1991discrete,rossler2009second}.}
    \label{tab:str_weak_accuracy_summary}
\end{table}

\section{Derivation and properties of benchmark problem}
\label{sec:benchmark}
To analytically compare the methods \eqref{eq:EM}, \eqref{eq:Mil}, \eqref{eq:Heu} and \eqref{eq:RK3}, we focus on a test problem where $f$ and $g$ are time-independent and affine in the spatial variable, i.e.
\[
f(\xi) = -\alpha -\beta \xi \quad\text{and}\quad g(\xi) = \gamma+\delta \xi.
\]
In other words, we consider the SDE
\begin{equation}\label{eq:GenericLinearSDE}
\rd \xi_\tau = -(\alpha +\beta \xi_\tau )\rd \tau+(\gamma+\delta \xi_\tau)\rd W_\tau ,
\end{equation}
which is the most general form of autonomous SDE in one variable with coefficients which are affine in $\xi_\tau$, driven by a single Brownian motion. This problem can be seen as the linearised form of the general nonlinear autonomous SDE \eqref{eq:Ito_SDE} about some fixed spatial point.

In this section, we perform what amounts to a non-dimensionalisation of this equation, which allows us to reduce the equation to a single parameter family in the most general case. We also discuss properties of the resulting reduced equation.

\subsection{Reduction of cases}\label{subsec:rescale}
We proceed by rescaling and shifting time and space coordinates from $(\xi, \tau) $ to $(x, t)$. To do so, we define $\xi_\tau:= Ax_t+B$, 
and we introduce the time rescaling $\tau = T t$. Here $A$, $B$ and $T$ are real parameters to be chosen, and $T>0$ is positive, so that the direction of time is preserved. Substituting appropriately, and using the properties of Brownian motion, we obtain the following equation in the new coordinates
\[
\rd x_t = -\left(T A^{-1}(\alpha+\beta B)+\beta T x_t \right) \rd t+T^{1/2}A^{-1}\left(\gamma+\delta B+\delta A x_t\right)\rd W_t.
\]
If $\beta\neq 0$, then we can set $B=-\alpha/\beta$ and $T=1/|\beta|$ and obtain
\[
\rd x_t = -\frac{\beta}{|\beta|}x_t\rd t+\left(\frac{\beta\gamma-\alpha\delta}{A\beta|\beta|^{1/2}} +\frac{\delta}{|\beta|^{1/2}}x_t\right)\rd W_t.
\]
Next, as long as $\beta\gamma-\alpha\delta\neq0$, we can set $A = \frac{\beta\gamma-\alpha\delta}{\beta|\beta|^{1/2}}$ and $\eta := \frac{\delta}{|\beta|^{1/2}}$ to obtain
\begin{equation}\label{eq:reducedSDE}
   \rd x_t = \pm x_t\rd t+(1+\eta x_t\big)\rd W_t,
\end{equation}
where the sign of the term in front of the drift term is opposite to that of $\beta$. Equation \eqref{eq:reducedSDE} provides the most generic dimensionless form of \eqref{eq:GenericLinearSDE} in this sense, but there are other possible cases:
\begin{itemize}
    \item If $\beta\neq0$ but $\beta\gamma-\alpha\delta = 0$, then we can rescale to obtain geometric Brownian motion
    \begin{equation}\label{eq:geomBM}
        \rd x_t = \pm x_t \rd t +\eta x_t \rd W_t,
    \end{equation}
    where again $\eta = \frac{\delta}{|\beta|^{1/2}}$.
    \item If $\beta=0$ but $\delta\neq0$, then we can set $B=-\gamma/\delta$, $T = |\delta|^{-1/2}$ and $A=-\alpha|\delta|^{1/2}$ to obtain
    \begin{equation}\label{eq:beta=zeroSDE}
        \rd x_t = \rd t\pm x_t\rd W_t.    
    \end{equation}
    \item Finally, we could have $\beta=0$ and $\delta=0$, resulting in an equation which is simply a translation of a Brownian motion.
\end{itemize}
In the existing literature on the numerical analysis of SDEs with multiplicative noise, geometric Brownian motion \eqref{eq:geomBM} has been used as a common test equation for numerical methods. A likely reason for the use of geometric Brownian motion is the availability of an analytical solution; the derivation above shows however that both this case and the case where \eqref{eq:beta=zeroSDE} is the resulting equation form sets of zero measure in the four-dimensional space of parameters $\alpha$, $\beta$, $\gamma$ and $\delta$. As such, we choose to focus on consideration of \eqref{eq:reducedSDE} for $\eta\in\R$ throughout the remainder of this work. We choose also to only consider the case with drift term being
$-x_{t}\rd t$ in \eqref{eq:reducedSDE}, which guarantees ergodicity of the dynamics with respect to an equilibrium distribution; this is discussed further in the following section.

\begin{rmrk}
As discussed in the introduction, we highlight some complexities which would arise in the more general case where we consider an SDE in $d$ dimensions, which motivate our choice to focus on a one-dimensional case. Consider the general $d$-dimensional system of SDEs where the drift and diffusion are affine, i.e.
\begin{equation*}
    \rd \xi^i_{\tau} = \left(\alpha^i+\sum_{j=1}^d\beta^{ij} \xi^i_{\tau}\right)\rd \tau + \sum_{j=1}^d\left(\gamma^{ij}+\sum_{k=1}^d\delta^{ijk} {\xi}^k_{\tau}\right)\rd W^j_{\tau},
\end{equation*}
where $\alpha^i$, $\beta^{ij}$, $\gamma^{ij}$ and $\delta^{ijk}$ are real coefficients of the system for $i,j,k\in\{1,\ldots,d\}$. In this case, the number of independent coefficients is $d+2d^2+d^3$. In higher dimensions, a general redefinition of the spatial variables takes the form
\begin{equation*}
   \xi^i_\tau := A^i+\sum_{j=1}^d B^{ij}{x}^j_{t},\quad i = 1,\dots, d,
\end{equation*}
along with the time rescaling $\tau = Tt$.
This provides us with the flexibility to scale out $d^2+d+1$ parameters in the most generic case, but under such a rescaling, there would remain $d^3+d^2-1$ independent dimensionless parameters governing the system. Even in two dimensions, there are therefore $11$ possible cases (although we note further reductions may be possible through finding other equivalences). This rapid growth in the number of cases with dimension justifies our choice to focus on a one-dimensional cases where there is exactly $1$ dimensionless variable for simplicity. Note that in the case where $d=1$, the above argument predicts that we have exactly one free dimensionless parameter, just as derived above.

For comparison, the particular case considered in \cite{saito2002mean} corresponds to a reduction of the two-dimensional case with the further requirement that $\beta$ is symmetric, $\gamma^{ij}$ is zero, $\delta^{ijk}=0$ for $k=2$, and $i,j\in\{1,2\}$. As in the case of geometric Brownian motion in one dimension, this is clearly not a case that will arise generically from a linearisation of nonlinear drift and diffusion about a fixed spatial point, but the analysis remains challenging: there are $9$ parameters that can be reduced to $2$ by scaling in space and time. This conclusion can effectively be seen from the result of Theorem~1 in \cite{saito2002mean}, which relies upon two independent conditions holding in order to guarantee stability.
\end{rmrk}

\subsection{Equilibrium distribution}
\label{sec:equilibrium}
Since our focus is on asymptotic statistical properties of the dynamics and the accurate recovery of these statistics by numerical schemes, we now consider the ergodicity properties of our reduced equation \eqref{eq:reducedSDE}. There are a range of established approaches which can be used to demonstrate that the SDE is ergodic with respect to an equilibrium distribution. One such approach is provided by the results of \cite{mattingly2002ergodicity}, which in turn relies upon the methodology introduced in \cite{down1995exponential,lord1995analysis,meyn2012markov}.
The Fokker--Planck equation which governs the probability density function $p$ of solutions to the model problem \eqref{eq:reducedSDE} is the PDE
\begin{equation}\label{eq:FP}
\partial_t p 
= \partial_x\Big(x p +\partial_x\left(\tfrac12(1+\eta x)^2 p\right)\Big).
\end{equation}
 In particular, this approach to showing ergodicity requires verification a \emph{minorisation condition} for the transition kernel of the process, and the construction of a \emph{Lyapunov function}. {In the case of the equation under consideration, the minorisation condition can be shown to hold as a consequence of H\"{o}rmander’s theorem \cite{H67,lord1995analysis}, and $V(x) = \sqrt{1+x^2}$ serves as an appropriate Lyapunov function.} However, a detailed verification of all conditions required to establish ergodicity is beyond the scope of the present work.

Instead, we focus here on formally computing the equilibrium distribution for use as a comparator for asymptotic statistics. When the SDE is ergodic, there exists a steady-state solution 
$p_\infty$ of \eqref{eq:FP} which is Lebesgue integrable on the real line such that, for most sufficiently regular initial distributions, the probability density $p_t$ converges to $p_{\infty}$ as $t\to +\infty$. Because $p_\infty$ is time independent, the steady state of \eqref{eq:FP} becomes 
\[
0 = \partial_x\Big(x p_\infty +\partial_x\left(\tfrac12(1+\eta x)^2p_\infty\right)\Big)
=\partial_x\Big(\left(\eta+(\eta^2+1) x\right)p_\infty+\tfrac12(1+\eta x)^2\partial_x p_\infty\Big),
\]
which holds wherever $p_\infty$ is supported.
We note that the coefficient of the highest order terms in this equation vanish when $x=-\tfrac1\eta$, and this reflects the compact support of the solution $p_\infty$ when $\eta\neq 0$.
Integrating once and rearranging, we have
\[
\frac{K}{(1+\eta x)^2} =\frac{2\eta+2(\eta^2+1)x}{(1+\eta x)^2}p_\infty+\partial_x p_\infty,
\]
where $K$ is some constant.
Multiplying through by the integrating factor
\[
I(x):=\exp\left(\int^x \frac{2\eta+2(\eta^2+ 1)u}{(1+\eta u)^2}du\right) = \exp\left(\frac{2}{\eta^2}\frac{1}{1+\eta x}\right)\Big(1+\eta x\Big)^{2+\frac{2}{\eta^2}},
\]
we have that
\[
K\exp\left(\frac{2}{\eta^2}\frac{1}{1+\eta x}\right)\Big(1+\eta x\Big)^{\frac{2}{\eta^2}}
=\partial_x\left(I(x)p_\infty(x)\right).
\]
We note that the term on the right hand side is not integrable in a neighborhood of $x=-\frac1\eta$, so we must have that $K=0$. Integrating again it follows that
\begin{equation*}
p_\infty(x) = \frac{C}{I(x)} = C\exp\left(-\frac{2}{\eta^2}\frac{1}{1+\eta x}\right)\Big(1+\eta x\Big)^{-2-\frac{2}{\eta^2}},
\end{equation*}
for some coefficient $C$, and to normalise appropriately, we can integrate. For now, assuming that $\eta>0$ and making the change of variable $y = \frac{2}{\eta^2}\frac{1}{1+\eta x}$, we can express the integral using the Gamma function as
\[
\int_{-\frac1\eta}^\infty \exp\left(-\frac{2}{\eta^2}\frac{1}{1+\eta x}\right)\Big(1+\eta x\Big)^{-2-\frac{2}{\eta^2}}\rd x = 2^{-\frac{2}{\eta^2}} \eta^{-1 + \frac{4}{\eta^2}} \Gamma(\tfrac{2}{\eta^2}).
\]
The density of the properly normalised form of the equilibrium distribution is therefore
\begin{equation}\label{eq:equilibrium}
p_\infty(x) =
\begin{cases}\displaystyle
    \frac{2^{\frac{2}{\eta^2}}\eta^{1-\frac{4}{\eta^2}}}{\Gamma(\frac{2}{\eta^2})}\exp\left(-\frac{2}{\eta^2}\frac{1}{1+\eta x}\right)\Big(1+\eta x\Big)^{-2-\frac{2}{\eta^2}} & x>-\frac{1}{\eta}, \\
    0 &x\leq -\frac1\eta.
\end{cases}
\end{equation}
If $\eta<0$, then we obtain the same form of the equilibrium distribution, but $x$ and $\eta$ are replaced by $-x$ and $-\eta$ in the definition \eqref{eq:equilibrium} above. In the special case where $\eta=0$, the SDE \eqref{eq:reducedSDE} has an additive noise term, and the distribution becomes a standard Gaussian with zero mean and variance $\frac{1}{2}$.

\begin{rmrk}\label{rmk:equilibrium_eta}
{We note that the resulting distribution \eqref{eq:equilibrium} is in fact a translated form of the inverse Gamma distribution, which is commonly parametrised with probability distribution function
\begin{equation*}
\tilde{p}(y;\alpha,\beta) = \frac{\beta^{\alpha}}{\Gamma(\alpha)}\frac{1}{x^{\alpha+1}}\exp(-\beta/y),
\end{equation*}
where $\alpha = 1+\frac{2}{\eta^2}$ and $\beta = \frac{1}{\eta^2}$, and $y$ is related to $x$ through the transformation $y=1+\eta x$.}
Due to the algebraic decay of the tail of this distribution as $x\to+\infty$, we note that the $k$th moment of the equilibrium distribution is only well-defined when $x^kp_\infty(x)$ is Lebesgue integrable, which requires
\[
k-2-\frac{2}{\eta^2}<-1,\quad\text{i.e.}\quad |\eta|<\sqrt{\frac{2}{k-1}}.
\]
{When they exist, we} can also explicitly calculate that the first and second moments of the equilibrium distribution. These are
\begin{equation}\label{eq:equilibrium_moments}
    \mu^{(1)}_\infty = \int_{-\infty}^\infty xp_\infty\,\rd x = 0\qquad\text{and}\qquad \mu^{(2)}_\infty = \int_{-\infty}^\infty  x^2p_\infty\,\rd x = \frac{1}{2-\eta^2};
\end{equation}
explicit representations as rational functions of $\eta$ can also be computed, but our focus will be on the first and second moments. We note that, when the second moment is finite, we can infer the value of $\eta$ up to its sign, and hence identify the equation governing the dynamics within the class of reduced equations \eqref{eq:reducedSDE}.
\end{rmrk}

\subsection{Evolution of moments}
\label{sec:moments}
In addition to finding the equilibrium distribution, our benchmark equation is simple enough that we are also able to use the Fokker--Planck equation \eqref{eq:FP} to find closed form expressions for the evolution of the moments of the distribution of independent solutions to \eqref{eq:reducedSDE} over time.
In particular, consider the evolution of the first and second moments of the position distribution, defining
\[
\mu^{(1)}(t) := \int_{-\infty}^\infty xp\,\rd x\quad\text{and}\quad \mu^{(2)}(t) := \int_{-\infty}^\infty x^2p\,\rd x.
\]
We can now use the Fokker--Planck equation \eqref{eq:FP}
to derive equations for these moments. For example, we have
\begin{equation*}
    \frac{\rd\mu^{(1)}}{\rd t} = \int_{-\infty}^\infty x\partial_tp\,\rd x = \int_{-\infty}^\infty x\partial_x\Big(x p +\partial_x\left(\tfrac12(1+\eta x)^2 p\right)\Big)\,\rd x.
\end{equation*}
Formally integrating by parts, we find that the first moment satisfies the ODE
\begin{equation*}
    \frac{\rd\mu^{(1)}}{\rd t} = -\int_{-\infty}^\infty x p \,\rd x -\int_{-\infty}^\infty \partial_x\left(\tfrac12(1+\eta x)^2 p\right)\rd x= -\mu^{(1)},
\end{equation*}
and a similar argument allows us to show that the second moment satisfies
\[
\frac{\rd\mu^{(2)}}{\rd t} = -(2-\eta^2)\mu^{(2)}+2\eta \mu^{(1)}+1.
\]
If the initial condition is deterministic, so that $p(x,0) = \delta_{x_0}(x)$, the corresponding initial conditions for these moment equations are $\mu^{(1)}(0)=x_0$ and $\mu^{(2)}(0) = x_0^2$. Solving the first equation, we find that
\begin{equation}\label{eq:analytic1stmomentevolution}
\mu^{(1)}(t) =x_0 \mathrm{e}^{-t}.
\end{equation}
As expected, $\mu^{(1)}(t)\to 0$ as $t\to\infty$, and the first moment tends to the first moment of the equilibrium distribution, as expected from the discussion of ergodicity given above.

For the second moment, we can substitute the expression for the first moment and solve, yielding the solutions
\begin{equation}
\mu^{(2)}(t) = \begin{cases}
    x_0^2\mathrm{e}^{-(2-\eta^2)t}+\frac{1}{2-\eta^2}(1-\mathrm{e}^{-(2-\eta^2)t})+\frac{2\eta x_0}{1-\eta^2} (\mathrm{e}^{-t}-\mathrm{e}^{-(2-\eta^2)t}), & \eta \neq \pm 1,\\
    x_0^2\mathrm{e}^{-t}+1-\mathrm{e}^{-t}\pm 2 x_0 t\mathrm{e}^{-t}, & \eta = \pm 1.
\end{cases}
\label{eq:analytic2ndmomentevolution}
\end{equation}
We observe that the long-time limit of $\mu^{(2)}$ exists only in the case when $\eta^2<2$, which corresponds to the case where the second moment of the equilibrium distribution is finite, and in this case $\mu^{(2)}(t)\to\frac{1}{2-\eta^2}$ as $t\to\infty$, which is exactly the value given in \eqref{eq:equilibrium_moments}. 

\begin{rmrk}\label{rmk:}
Even in the case where $\mu^{(2)}$ blows up as $t\to\infty$, we observe that we can still determine the value of $\eta$ (up to sign) by observing the exponential growth rate:
\[
\lim_{t\to\infty}\frac{\log \mu^{(2)}(t)}{t} = \eta^2-2.
\]
This demonstrates that first and second moment information are sufficient to recover the $\eta$ parameter which fixes the particular reduced problem \eqref{eq:reducedSDE}.
\end{rmrk}

\section{Statistical stability properties of numerical schemes}
\label{sec:stab}
We now turn our focus to the asymptotic stability properties of the numerical schemes we consider.
To study the asymptotic properties of the numerical schemes considered across a range of time-step sizes $h$, we define the $j$th moment at step $n$ of the numerical scheme to be
\[
\mu_n^{(j)}:=\mathbb{E}[(x_n)^j],
\]
where $x_n$ is the position predicted by one of the schemes described in Section~\ref{sec:schemes} at time $t_n = nh$.
{
By considering each of the schemes in turn, we will show that we can derive recurrence relations for these moments. We will focus on the cases where $j=1$ and $j=2$, since as we have argued above, these are the moments which together can be used to characterise the reduced SDE \eqref{eq:reducedSDE} via the method of moments. Moreover, in a general nonlinear setting, it is common to focus on wide-sense (or weakly) stationary processes; these are the class of processes which have time-invariant mean and time autocovariance functions, rather than being strictly stationary (in which case all finite time-marginals are time-invariant). Wide-sense stationary processes are often considered since the mean, variance and autocovariance are the most important and readily-measurable statistics for applications, and being wide-sense stationary is equivalent to being strictly stationary for Gaussian processes. The latter are a common model for many data-generating processes in practice; see for example the discussion in Section~4.18 of \cite{GD04}.}

\subsection{Asymptotic moment stability}
For the purposes of our analysis, we define the following notions of asymptotic stability and accuracy of a numerical discretisation as applied to a particle SDE problem.

\begin{dfntn}
    For a given time step $h$, a numerical discretisation scheme of an SDE is said to be \emph{asymptotically stable for the $j$th moment} if
        \[
            \limsup_{n\to \infty}\left|\mu^{(j)}_n-\mathbb{E}[(X_{t_n})^j]\right|<+\infty.
        \]
    A numerical discretisation is said to be \emph{asymptotically $O(h^n)$ accurate for the $j$th moment} if for all $h$ sufficiently small,
    \[
        \limsup_{n\to\infty}\left|\mu^{(j)}_n-\mathbb{E}[(X_{t_n})^j]\right| = O(h^n).
    \]
\end{dfntn}

\noindent
Clearly, the latter notion of asymptotic accuracy for the $j$th moment necessarily requires that there exist some possible values of the time step $h$ for which the scheme is also asymptotically stable for the same moment.

{
Before proceeding further, we comment on the distinction between the definition of asymptotic moment stability given above and the existing condition of mean-squared (MS) stability, which was first defined in \cite{saito1996stability}. In fact, our definition is simply a generalisation of the notion of MS-stability: MS-stability was introduced on the basis of geometric Brownian motion as a benchmark, and corresponds to the requirement that $\mu^{(2)}_n\to 0$ as $n\to\infty$ for fixed step size $h$. The choice of this condition was motivated by the fact that it is convenient to understand a test problem in which all sample paths concentrate at $0$ as $t\to\infty$, as they do in the case of the benchmark ordinary differential equation $\rd x_t = -x_t\rd t$. In the case of geometric Brownian motion as studied in \cite{saito2002mean}, we have that
$\lim_{n\to\infty} \mathbb{E}[(X_{t_n})^2] = 0$. In this case, under our definition, an integrator applied to the same equations with fixed timestep $h$ is both stable and $o(1)$ accurate for second moments if and only if it is MS-stable.

By refining our understanding of stability and accuracy via the definition above, we can consider the discretisation of SDEs where $j$th moments tend to a stable non-zero value as $t\to\infty$. As explained in the introduction, this is often a desirable property of SDEs used in applications, where some form of ergodicity is of interest. Moreover, as we have already seen in the discussion of our benchmark given in Section~\ref{sec:equilibrium}, there are examples of relatively simple SDEs in which $\mathbb{E}[(X_{t_n})^j]$ does not have a finite limiting value, despite the SDE being ergodic. Our definition also covers such cases, by requiring only that the limit supremum of the difference between moments computed from the integrator and from the true SDE remain under control.}

{Over the course of the following subsections, we derive conditions under which the schemes described in Section~\ref{sec:schemes} are asymptotically stable for first and second moments when applied to the model problem \eqref{eq:reducedSDE}. A corresponding analysis for the more standard GBM benchmark can be found in Appendix~\ref{app:GBM}.}

\subsection{Euler-Maruyama (EM)}
\label{sec:EMstab}
{We begin by deriving a recurrence relation for the moments as they evolve under the EM scheme.} Taking expectations and using the scheme definition \eqref{eq:EM} in this case, we find that
\begin{equation*}
    \mu^{(1)}_{n+1}=\mathbb{E}[x_{n+1}] = \mathbb{E}\left[(1-h)x_n+(1+\eta x_n)\Delta W_n\right]= (1-h)\mu^{(1)}_n,
\end{equation*}
where we have used the fact that $x_n$ and $\Delta W_n\sim \mathcal{N}(0,\,h)$ are independent, and hence we can solve to obtain $\mu^{(1)}_n = (1-h)^n x_0$.
For the expectation to vanish in the long-time limit, we therefore require that
\[
|1- h|<1,\quad\text{i.e.}\quad 0<h<2.
\]
As we should expect, this is the classical stability region for the explicit Euler method.
When the scheme is asymptotically first moment stable, we have that $\mu^{(1)}_n\to0$ as $n\to\infty$, so there is no asymptotic bias in the scheme.

Following a similar approach for the second moment using independence and the properties of a standard normal random variable to deduce that the cross terms vanish, we have
\begin{align*}
\mu^{(2)}_{n+1}=\mathbb{E}[(x_{n+1})^2] &= \mathbb{E}[(1- h)^2x_n^2]+\mathbb{E}[(1+\eta x_n)^2(\Delta W_n)^2]\\
&=(1-h)^2 \mu^{(2)}_n+h(1+2\eta \mu^{(1)}_n+\eta^2\mu^{(2)}_n)\\
&= \left((1- h)^2+h\eta^2\right)\mu^{(2)}_n +2h\eta \mu^{(1)}_n+ h.
\end{align*}
Employing the standard linear stability analysis for discrete dynamical systems, the fixed point is stable when
\[
\left|\frac{\partial}{\partial \mu^{(2)}_n}\left[\left((1- h)^2+h\eta^2\right)\mu^{(2)}_n+2h\eta \mu^{(1)}_n+ h\right]\right| = \left|(1- h)^2+h\eta^2\right|<1;
\]
this reduces to the requirement that
\begin{equation}\label{eq:EM_mu2_stab}
    0<h<2-\eta^2.
\end{equation}

Note that this is (as expected) strictly more restrictive than the requirement that $0<h<2$ when only considering the asymptotic stability of the first moment.

If the scheme is asymptotically first and second moment stable under the conditions obtained above, then $\mu^{(2)}_n\to\mu^{(2)}_\infty$ which solves
\[
\mu^{(2)}_\infty = \left((1- h)^2+h\eta^2\right)\mu^{(2)}_\infty+ h\quad\text{and hence}\quad
\mu^{(2)}_\infty =\frac{1}{2-\eta^2-h} = \E[X_{t_n}^2]+O(h).
\]
The asymptotic second moment is therefore asymptotically first-order accurate for the second moment.

\subsection{Milstein (Mil)} For the Milstein scheme applied to test equation \eqref{eq:reducedSDE}, the relevant coefficients are
\[
f_n = -x_n,\quad g_n = 1+\eta x_n\quad\text{and}\quad \tfrac12[g'g]_n = \tfrac12\eta(1+\eta x_n),
\]
so
\[
x_{n+1} = (1-h)x_n+(1+\eta x_n)\Delta W_n+\tfrac12\eta(1+\eta x_n)(\Delta W_n^2-h).
\]
Taking expectations, the latter terms on the right-hand side vanish, and as with the Euler-Maruyama scheme, we have
\[
\mu^{(1)}_{n+1} = (1- h) \mu^{(1)}_n.
\]
The same step-size restriction therefore applies for asymptotic stability of the first moment, i.e. we require $0<h<2$. 
Similarly, when asymptotically first moment stable, there is no error in the scheme; the limit point is $\mu^{(1)}_\infty=0$.

To analyse the second moment, we find that after squaring, the expectations of the cross terms again vanish, and we obtain
\begin{align*}
    \mu^{(2)}_{n+1} &= (1- h)^2\mathbb{E}[x_n^2]+\mathbb{E}[(1+\eta x_n)^2]\E[\Delta W_n^2]+
    \tfrac14\eta^2\mathbb{E}[(1+\eta x_n)^2]\mathbb{E}[(\Delta W_n^2-h)^2]\\
    &=(1-h)^2\mu^{(2)}_n+h\big(1+2\eta\mu^{(1)}_n
    +\eta^2\mu^{(2)}_n\big)
    +\tfrac12h^2\eta^2\big(1+2\eta\mu^{(1)}_n+\eta^2\mu^{(2)}_n\big).
\end{align*}
Once more, assuming that $\mu^{(1)}_n\to0$, the system has a stable fixed point under the condition that
\begin{multline*}
    \left|\frac{\partial}{\partial \mu^{(2)}_n}\left((1-h)^2\mu^{(2)}_n+h\big(1+2\eta\mu^{(1)}_n
    +\eta^2\mu^{(2)}_n\big)
    +\tfrac12h^2\eta^2\big(1+2\eta\mu^{(1)}_n+\eta^2\mu^{(2)}_n\big)\right)\right|\\
    =\left|(1-h)^2+h\eta^2+\tfrac12h^2\eta^4\right|=\left|1-(2-\eta^2)h+(1+\tfrac12\eta^4)h^2\right|<1.
\end{multline*}
For $h>0$, this condition can be further reduced to
\begin{equation}\label{eq:Mil_mu2_stab}
    0<h<\frac{2-\eta^2}{1+\tfrac12\eta^4}.    
\end{equation}
Note that the latter condition for the stability of the second moment is more restrictive than for EM; compare \eqref{eq:EM_mu2_stab}.
When the scheme is asymptotically second moment stable, the fixed point $\mu^{(2)}_\infty$ must solve
\[
\mu^{(2)}_\infty = \left((1-h)^2+h\eta^2+ \tfrac12 h^2\eta^4\right) \mu^{(2)}_\infty + h +\tfrac12h^2\eta^2,
\]
which has solution
\[
\mu^{(2)}_\infty = 
\frac{1+\tfrac12h\eta^2}{2-\eta^2-h(1+ \tfrac12\eta^4)} = \frac{1}{2-\eta^2}+O(h),
\]
and hence the scheme is asymptotically first-order accurate for the second moment.

\subsection{Stochastic Heun (SH)}
To analyse the SH method, we note that the auxiliary function used in the scheme $F:=f-\tfrac12g'g$ is
\[
F(x) = -x-\tfrac12\eta(1+\eta x).
\]
Working through the algebra, and taking expectations, we find that
\[
\mu^{(1)}_{n+1} =  \left(1-h+\tfrac18h^2(2+\eta^2)^2\right)\mu^{(1)}_n
+\tfrac18\eta h^2(2+\eta^2).
\]
In this case, the criterion for first moment stability is
\[
\left|1-h+\tfrac18h^2(2+\eta^2)^2\right|<1,\quad\text{or}\quad 0<h<\frac{8}{(2 + \eta^2)^2}.
\]
When stable, the limiting fixed point $\mu^{(1)}_\infty$ is
\[
\mu^{(1)}_\infty =\frac{\tfrac18\eta h(2+\eta^2)}{1-\tfrac18h(2+\eta^2)^2} = O(h).
\]
As such, the SH scheme is asymptotically biased, i.e. it has a first-order accurate asymptotic first moment.

For the second moment, working through the algebra, we obtain the recurrence relation:
\begin{align*}
    \mu^{(2)}_{n+1}&=
    \Big(
    h 
    -\tfrac{1}{2}h^2 (2 + \eta^2)
    +\tfrac{1}{4} h^3 (1 + \eta^2)^2
    +\tfrac{1}{64}h^4\eta^2 (2 + \eta^2)^2
    \Big)\\
    &\qquad
    +\left(
    2h\eta
    -\tfrac{1}{4}h^2\eta(10 + 3 \eta^2)
    +\tfrac{1}{4} h^3 \eta(2 + 5 \eta^2 + 2 \eta^4)
    +\tfrac{1}{32} h^4 \eta(2 + \eta^2)^3 \right)\mu^{(1)}_n
    \\
    &\qquad+
    \Big(
    1
    -h(2 - \eta^2) 
    -\tfrac{1}{4}h^2 \big((2+\eta^2)^2-12\big)\\
    &\qquad\qquad \qquad
    +\tfrac{1}{4}h^3 (\eta^6+3\eta^4-4)
    +\tfrac{1}{64}h^4 (2 + \eta^2)^4
    \Big)\mu^{(2)}_n.
\end{align*}
This recurrence relation is stable if the first moment converges, and if in addition
\[
\left|1
    -h(2 - \eta^2) 
    -\tfrac{1}{4}h^2 \big((2+\eta^2)^2-12\big)
    +\tfrac{1}{4}h^3 (\eta^6+3\eta^4-4)
    +\tfrac{1}{64}h^4 (2 + \eta^2)^4
    \right| < 1.
\]
In this case, the stability region does not have simple closed form expression, and indeed we find that the SH method is in fact stable beyond the region in which the SDE has an asymptotically stable second moment; see Section~\ref{sec:num_verif} below for further discussion.

Using the asymptotic first moment obtained above and solving to find the fixed point, we obtain that the second moment is first-order accurate, i.e. that
\[
\mu^{(2)}_\infty=\frac{1}{2-\eta^2}+O(h).
\]

\subsection{3-stage Runge-Kutta}
\label{sec:RK3stab}
Here, after working through the algebra and taking expectations, we find that
\[
\mu^{(1)}_{n+1}=
    -\tfrac{1}{24} h^2 \eta(2 + 3 \eta^2)
    -\tfrac{1}{48} h^3 \eta(2 + \eta^2)^2
    +
    \Big(
    1
    -h
    +\tfrac{1}{8} h^2 (4 - \eta^4)
    -\tfrac{1}{48} h^3 (2+\eta^2)^3
    \Big) \mu^{(1)}_n
\]
The first moment is asymptotically stable whenever
\[
    \Big|
    1
    -h
    +\tfrac{1}{8} h^2 (4 - \eta^4)
    -\tfrac{1}{48} h^3 (2+\eta^2)^3
    \Big|
    <1,
\]
and the fixed point is
\[
\mu^{(1)}_\infty =
    -\frac{h \eta (4 + 6 \eta^2 + h (2 + \eta^2)^2)}{
 48 - 6 h (4 - \eta^4) + h^2 (2 + \eta^2)^3}
 =O(h).
\]
We therefore see that, like the Heun method, the Runge-Kutta scheme is asymptotically biased, but is asymptotically first-order accurate for the first moment.

For the second moment, we have
\begin{align*}
    \mu^{(2)}_{n+1}&=
    \Big(
    h
    -\tfrac{1}{2} h^2 (2 - \eta^2)
    +\tfrac{1}{12}h^3 (8 - \eta^4) 
    -\tfrac{1}{192}h^4 (32 + 20 \eta^2 - 44 \eta^4 - 27 \eta^6)\\
    &\qquad
    +\tfrac{1}{288}h^5 (2 + \eta^2)^2 (2 + 7 \eta^2 + 6 \eta^4)
    +\tfrac{1}{2304}h^6 \eta^2 (2 + \eta^2)^4
    \Big)\\
    &\quad+
  \Big(
    2 h \eta 
    -\tfrac{1}{12} h^2 \eta  (38 - 9\eta^2)
    +\tfrac{1}{24} h^3 \eta  (56 + 2 \eta^2 - 5\eta^4)\\
    &\qquad \quad 
    -\tfrac{1}{96} h^4 \eta  (72 + 44 \eta^2 - 50\eta^4 - 27\eta^6)\\
    & \qquad \qquad \qquad
    +\tfrac{1}{72} h^5 \eta (2 + \eta^2)^3 (1 + 3 \eta^2)
    +\tfrac{1}{1152} h^6 \eta (2 + \eta^2)^5
    \Big) \mu^{(1)}_n\\
    &\quad+\Big(
    1
    -h (2-\eta^2)
    +\tfrac{1}{4} h^2 (8 - 8\eta^2 + \eta^4)
    -\tfrac{1}{24} h^3 (32 - 36\eta^2 + 3\eta^6)\\
    &\qquad\quad
    +\tfrac{1}{192} h^4 (2 + \eta^2)^2 (28 - 52\eta^2 + 27\eta^4)
    -\tfrac{1}{96} h^5 (2 + \eta^2)^4 (1 - 2\eta^2)\\
    &\qquad \qquad \qquad 
    +\tfrac{1}{2304} h^6 (2 + \eta^2)^6
    \Big) \mu^{(2)}_n.
\end{align*}
The stability criterion in this case is
\begin{multline*}
    \Big|
        1
    -h (2-\eta^2)
    +\tfrac{1}{4} h^2 (8 - 8\eta^2 + \eta^4)
    -\tfrac{1}{24} h^3 (32 - 36\eta^2 + 3\eta^6)\\
    +\tfrac{1}{192} h^4 (2 + \eta^2)^2 (28 - 52\eta^2 + 27\eta^4)
    -\tfrac{1}{96} h^5 (2 + \eta^2)^4 (1 - 2\eta^2)
    +\tfrac{1}{2304} h^6 (2 + \eta^2)^6
    \Big|<1.
\end{multline*}
Again, there is no simple closed form for the solution in this case.

Under the assumption that the first moment is stable, it can be checked that the fixed point for this recurrence relation is
\[
\mu^{(2)}_\infty = \frac{1}{2-\eta^2}+O(h);
\]
the full expression is omitted due to its length, but we see that this method is asymptotically first-order accurate for the second moment.

\subsection{Summary of results}

\begin{table}[tp!]
    \centering
    \begin{tabular}{r|c|c}
        \textbf{Method} & \textbf{1st moment} & \textbf{2nd moment} \\
        \hline
        Euler-Maruyama (EM) &  Exact & $O(h)$\\
        Milstein (Mil) & Exact & $O(h)$ \\
        Stochastic Heun (SH) & $O(h)$ & $O(h)$ \\
        3-stage Runge-Kutta (RK3) & $O(h)$ & $O(h)$
    \end{tabular}
    \caption{A summary of the results of in Sections~\ref{sec:EMstab}--\ref{sec:RK3stab}, i.e. the asymptotic accuracy results for the first and second moments as predicted by the methods studied.}
    \label{tab:accuracy_summary}
\end{table}

{The asymptotic statistical accuracy results we have obtained over the course of the previous sections are summarised in Table~\ref{tab:accuracy_summary}, and a comparison of the stability regions for the first and second moments is provided in Figure~\ref{Fig:stab}.

Concerning the asymptotic accuracy results in Table~\ref{tab:accuracy_summary}, we recall that the results of \cite{ruemelin1982numerical,saito1991discrete,kloeden2012numerical} provide the strong and weak order of accuracy of various schemes over a fixed time interval $[0,T]$, and were summarised in Table~\ref{tab:str_weak_accuracy_summary}. While one might assume \emph{a priori} that schemes with higher weak order should give improved convergence of moments since these are expressed as expectations, we see that there are marked discrepancies between the sets of results. These gaps arise due to the time-dependent prefactor which is typically included in convergence estimates when considering a scheme over a fixed time interval $[0,T]$. As $T\to\infty$, this prefactor can blow up, meaning that the convergence estimates break down, and hence we do not necessarily achieve the same weak accuracy on asymptotic moments the the results on a finite time interval would predict. As such, we have shown the surprising result that selecting a higher-order scheme does not guarantee better approximation regarding the accuracy of asymptotic moments in the long-time simulation, at least when measured purely in terms of the time-step $h$, and indeed, higher-order schemes can in fact be less numerically stable in recovering these moments in long-term simulations.

\begin{figure}[tp!]
\centering
\includegraphics[width = \linewidth]{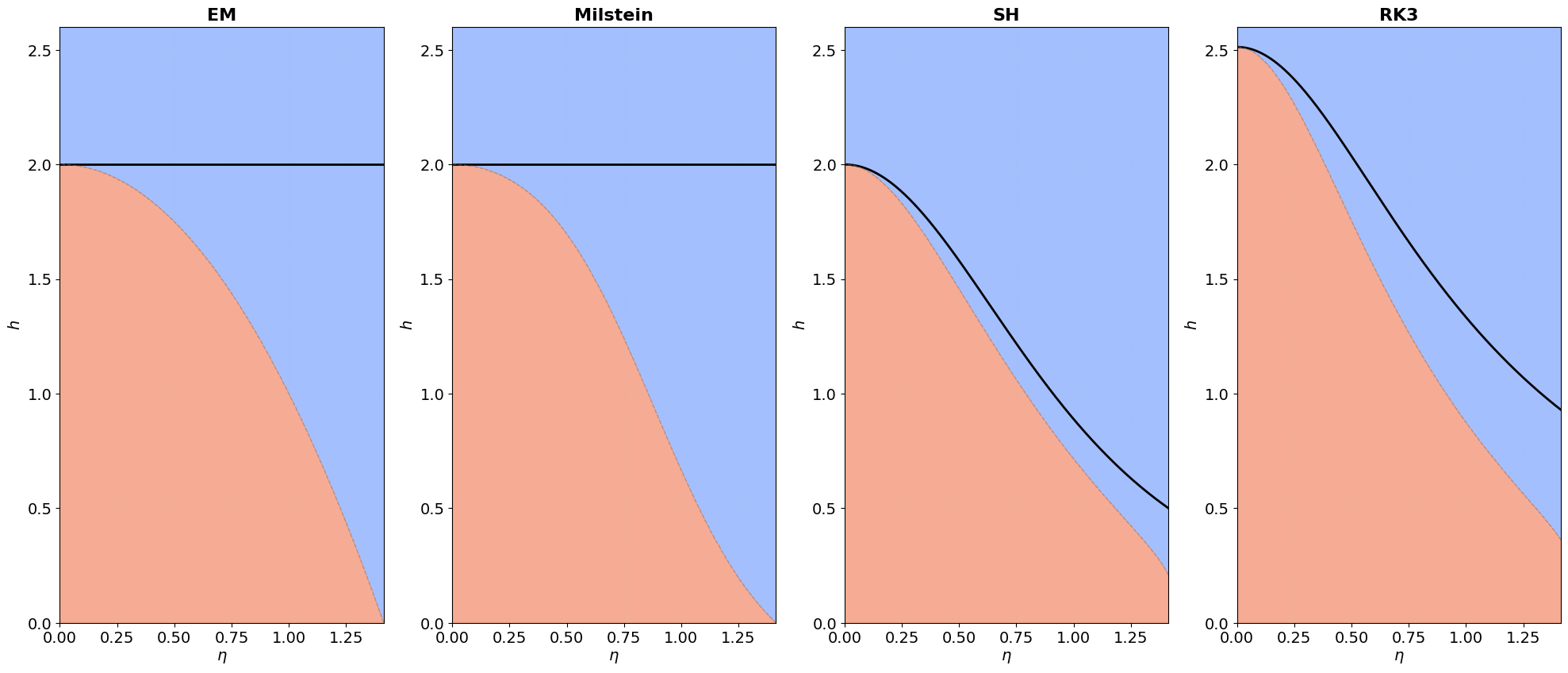}\\
\includegraphics[width = 0.5\linewidth]{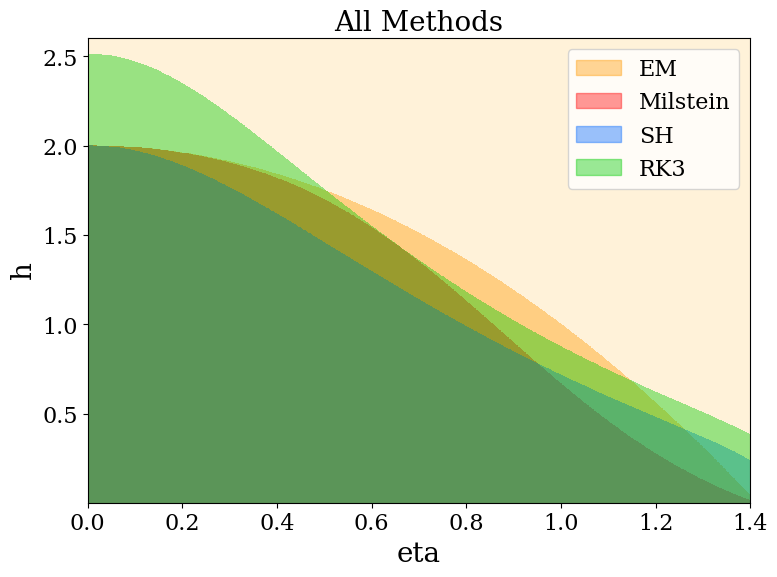}
\caption{\textbf{Top row}: 1st moment stability regions (boundaries shown as lines) and 2nd moment stability regions (pink shaded regions) for the schemes considered; in each case, the stability region for second moments is strictly contained inside the stability region for first moments. \textbf{Bottom row}: A comparison of 2nd moment stability regions among all schemes plotted on the same axes. For asymptotic first moments, we observe EM and Mil are more stable than the higher-order schemes in many cases. For asymptotic second moments, RK3 is the most stable scheme when $\eta\le 0.5$ or $\eta > 1.145$; and EM is the most stable when $0.5 < \eta < 1.145$.  }
\label{Fig:stab}
\end{figure}

Turning to the stability regions shown in Figure~\ref{Fig:stab}, we see that for first moments, EM and Mil have identical stability regions, both regions being larger than SH, whereas RK3 is more stable than EM and Mil for $\eta\lesssim 0.5245$, and less stable if $\eta\gtrsim 0.5245$. This starkly illustrates the fact that when $\eta$ is reasonably large, the lower-order schemes recover bounded first moments over a much wider range of time-steps $h$, while the higher-order schemes become unstable in long time simulations. As such, we see there is not always a stability gain to be made from using a higher-order method, in agreement with our analysis.

As our analysis predicts, we see that the regions in parameter space over which second moments are stable lie inside the first moment regions, and these regions coincide as $\eta\to0$, which is the case where the diffusion coefficient is no longer position dependent.
For 2nd moments, we observe that EM is stable over a strictly larger region than Mil, and RK3 is stable over a strictly larger region than SH. When $\eta$ is very small or very large, that is $\eta\le 0.5$ or $\eta\ge 1.145$, RK3 is more stable than EM, however, for $0.5<\eta<1.145$, EM is better. For second moments, we only see significant gains in the range of stable step sizes when $\eta$ approaches $\sqrt{2}$, but we note that the bias observed in the first moments for the higher-order schemes would lead to poor recovery of the variance even in this case.

Since the geometric Brownian motion (GBM) benchmark is more standard in the analysis of SDE schemes than our benchmark \eqref{eq:reducedSDE} and was used in the original definition of the notion of MS-stability \cite{saito1996stability}, in Appendix~\ref{app:GBM}, we present a parallel analysis for the application of the schemes to our rescaled form of GBM \eqref{eq:geomBM}. In terms of the timestep $h$ and dimensionless parameter $\eta$ present in this equation, we find that the stability regions take the same form, and so Figure~\ref{Fig:stab} also illustrates these regions. We note however that a stark difference between the analysis of schemes applied to our proposed benchmark \eqref{eq:reducedSDE} and GBM is that all schemes are asymptotically exact in terms of the first and second moment for GBM. This highlights a special feature of the GBM equation, which is that the drift and diffusion both vanish at the same point, and this means that, when it exists, the equilibrium measure is concentrated at $x=0$. As argued in our derivation in Section~\ref{subsec:rescale}, this behaviour is rather special when considering a general SDE with affine coefficients, and is not always typical in applications. Moreover, it hides the variability in long-term performance; this will become clearer when we consider a nonlinear example equation in Section~\ref{sec:nonlinear}.}

\section{Numerical verification}
\label{sec:num_verif}
To verify and illustrate our analytical results, we performed a range of numerical tests, both on the benchmark example equation, and on a nonlinear equation. These results are presented over the course of this section.

\subsection{Asymptotic accuracy of moments for methods considered}

\begin{figure}[tp!]
\centering 
\subfigure
{\includegraphics[width = 0.46\linewidth]{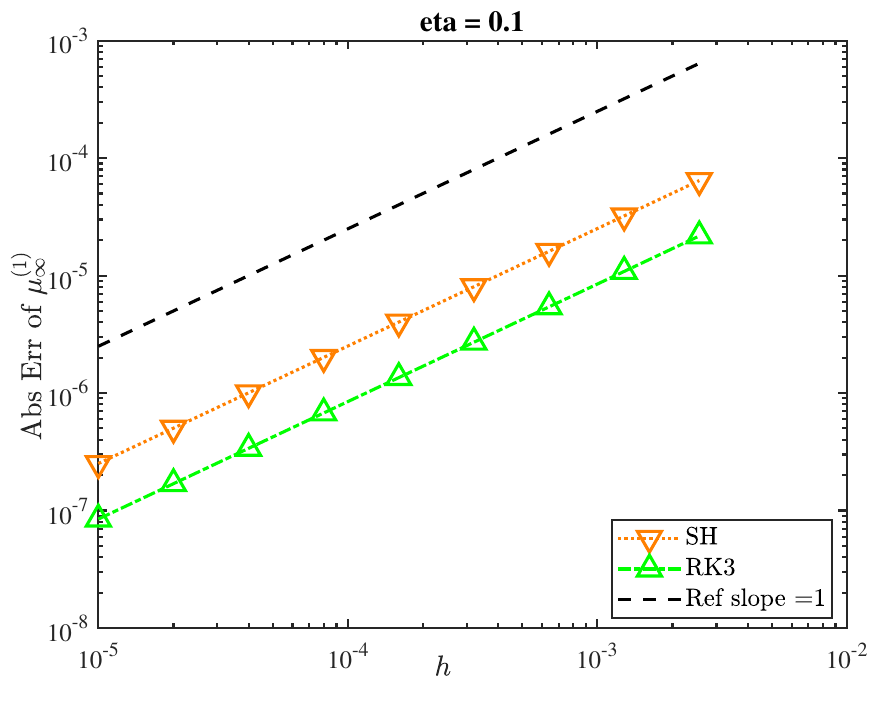}}
\subfigure
{\includegraphics[width = 0.46\linewidth]{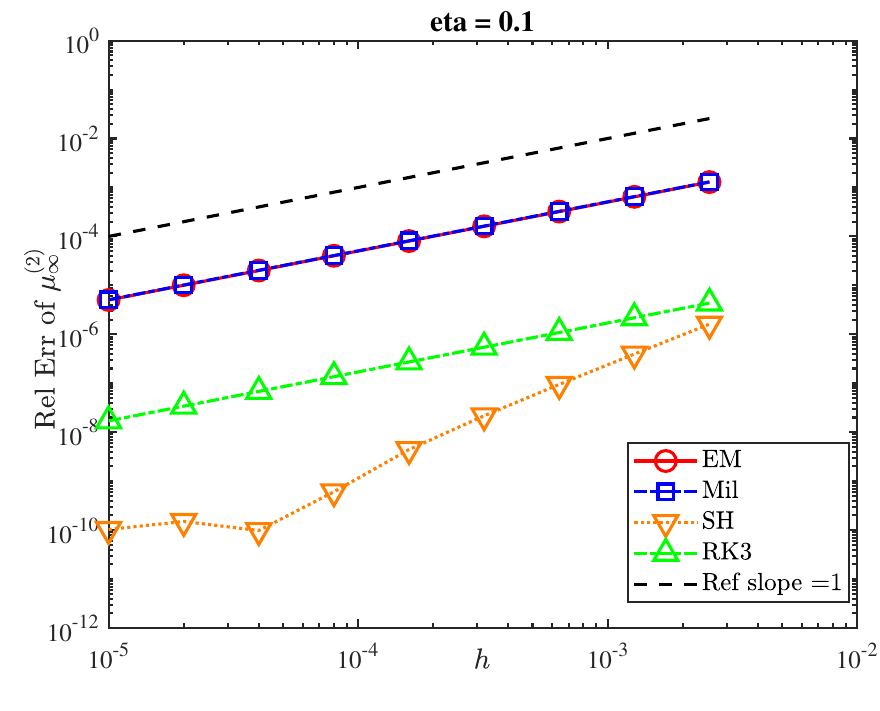}}
\subfigure
{\includegraphics[width = 0.46\linewidth]{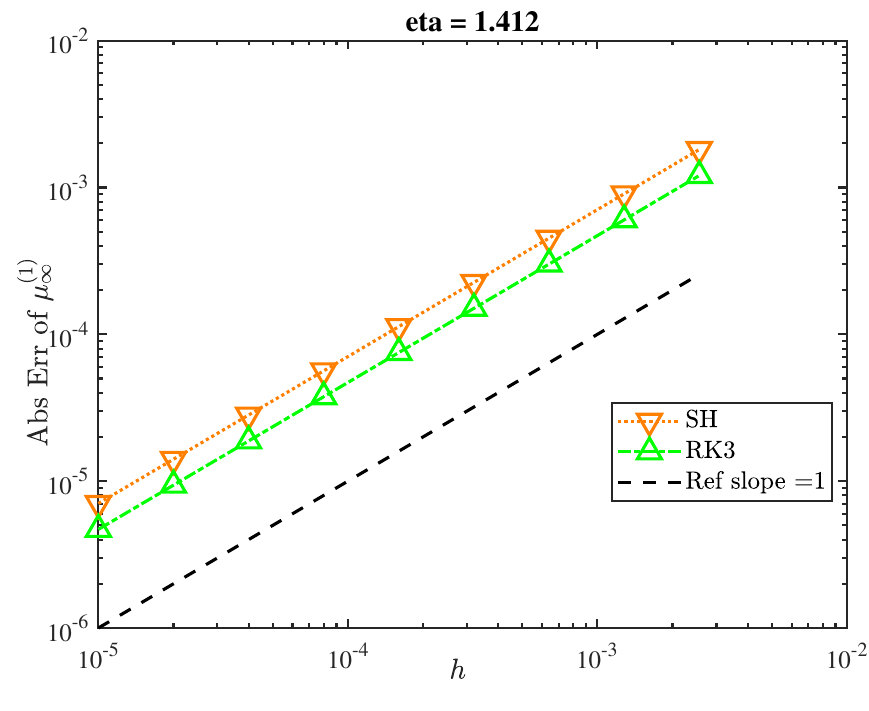}}
\subfigure
{\includegraphics[width = 0.46\linewidth]{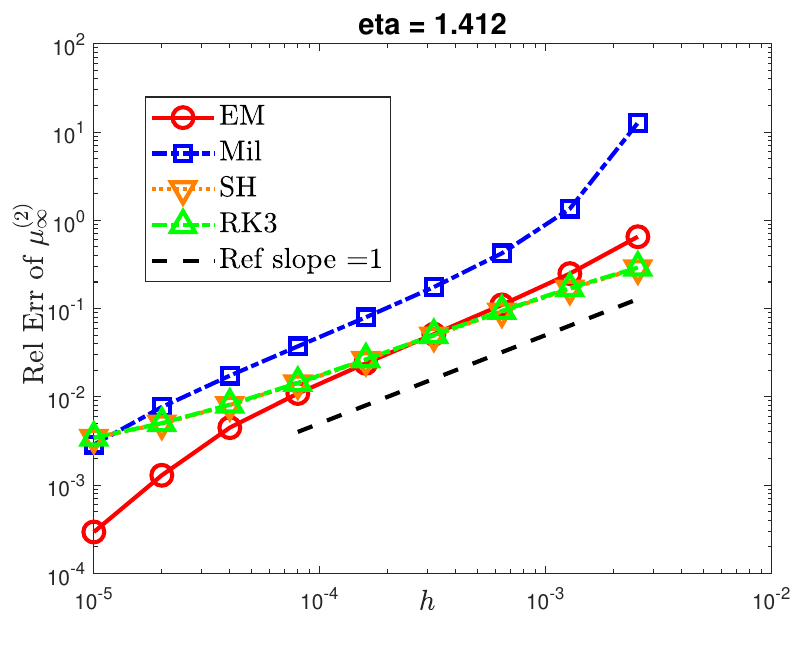}}
\caption{Asymptotic accuracy of the first and second moments obtained using the various numerical methods with respect to the step size $h$. Errors were computed as the absolute value of the difference between the exact moment in the limit $T\to\infty$ and the numerical moment at final simulation time $T=1000$.
\textbf{Column 1:} The absolute error in the first moment for the SH and RK3 schemes only, showing the $O(h)$ convergence predicted; EM and Mil are not shown as they exactly recover the first moment. \textbf{Column~2:} The relative error of the asymptotic second moment for each scheme. For the case where $\eta=0.1$ (first row), the higher-order schemes perform better, and SH reaches machine precision for small $h$, with the accumulation of floating point errors appearing to limit further gains. For the case of larger $\eta$ (second row), all schemes give similar error, with EM outperforming the higher-order schemes as $h$ gets closer to $0$.}\label{Fig:accuracy}
\end{figure}

To verify the asymptotic accuracy results summarized in Table~\ref{tab:accuracy_summary}, we computed the long-time errors of the first and second moments obtained by various numerical methods with respect to the step size $h$, and compared against the asymptotic moments of the equilibrium distribution: these results are shown in Figure~\ref{Fig:accuracy}.

Recall that the reference exact asymptotic moments are given by
\[
\mu_{\infty}^{(1)} =0\quad\text{and}\quad \mu_{\infty}^{(2)} =\frac{1}{2-\eta^2},
\]
so we set the final simulation time to $T=1000$ and considered two different parameter values for the dimensionless parameter $\eta$ where we expect different types of behaviour. The values chosen were $\eta = 0.1$ and $\eta=1.412$, respectively. We see that errors decrease consistently and approximately linearly in most cases as the time step $h$ decreases, following the estimates summarized in Table~\ref{tab:accuracy_summary}. For the small $\eta=0.1$, we note that the SH method appears to exhibit superlinear convergence rate for the second moment; we hypothesise that this may be due to improved pre-asymptotic performance of this scheme in the case where the diffusion term is close to being position independent. In contrast, for the larger value 
$\eta=1.412$, the EM scheme is very competitive with the higher-order schemes in recovering the second moment, but all schemes show similar performance at moderate time step sizes.

\subsection{Pre-asymptotic moment evolution}
\begin{figure}[htp!]
    \centering
   \includegraphics[width=0.92\linewidth]{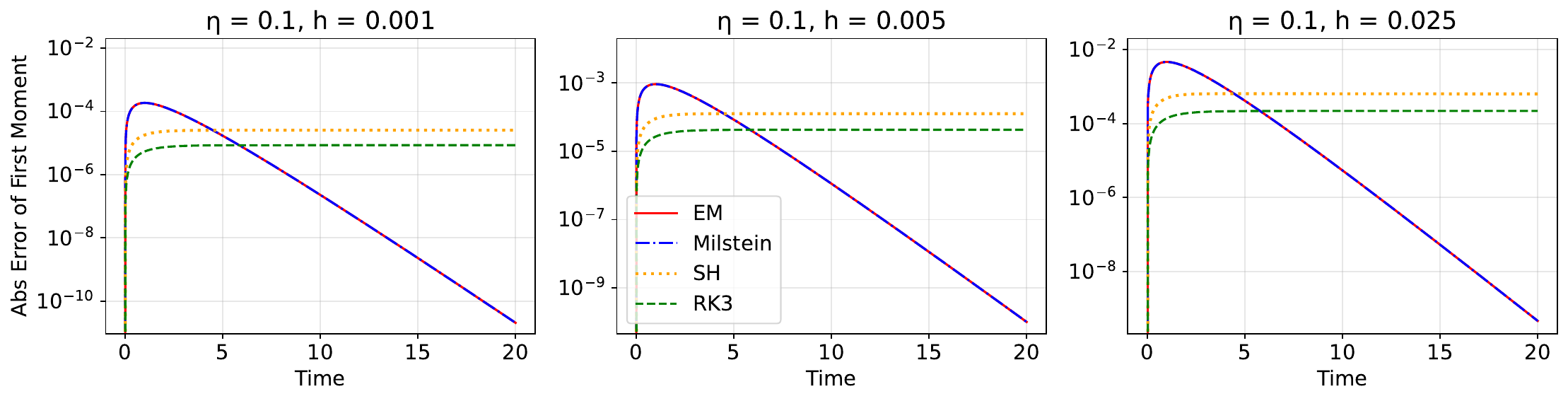}
   \includegraphics[width=0.92\linewidth]{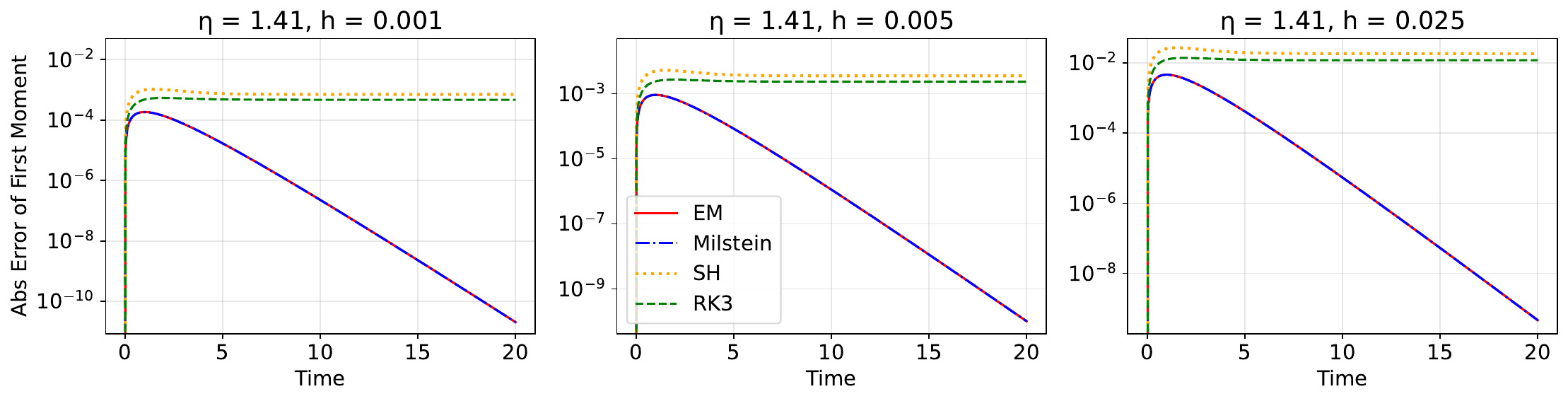}
    \caption{A comparison of the evolution of absolute first moment error for the various discretisation schemes up to time $T=20$ for our benchmark equation. The first row depicts the case when $\eta=0.1$ for step sizes $h=0.005,0.005,0.025$, and the second row depicts when $\eta=1.41$ with the same range of step sizes. In agreement with our analytical results, we see that the EM and Mil recover asymptotically exact first moments, whereas RK3 and SH are asymptotically biased.
    }
    \label{fig:1stMom}
\end{figure}

\begin{figure}[htp!]
    \centering
    \includegraphics[width=0.92\linewidth]{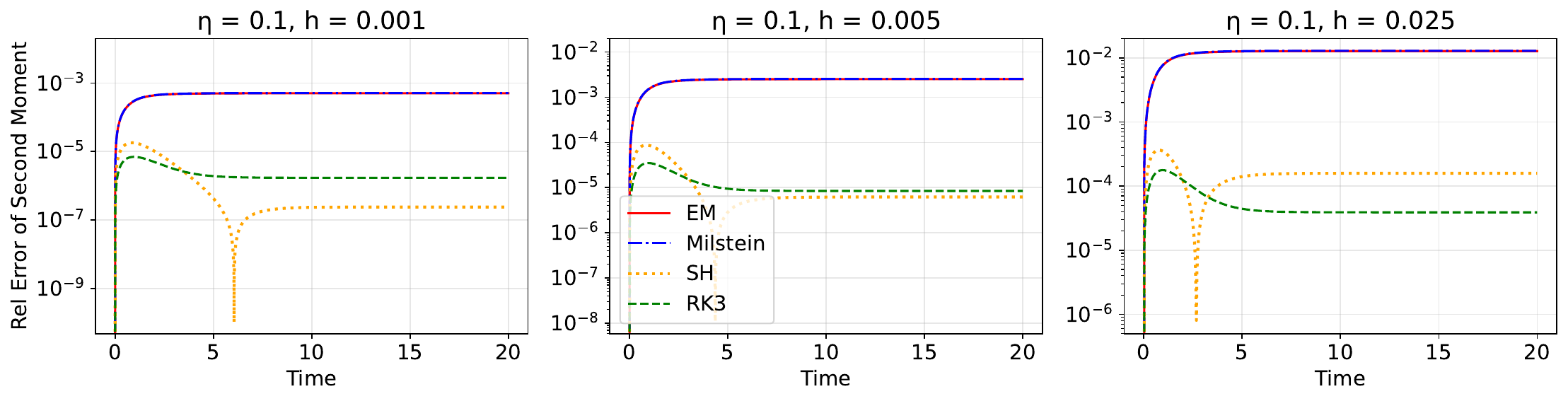}
    \includegraphics[width=0.92\linewidth]{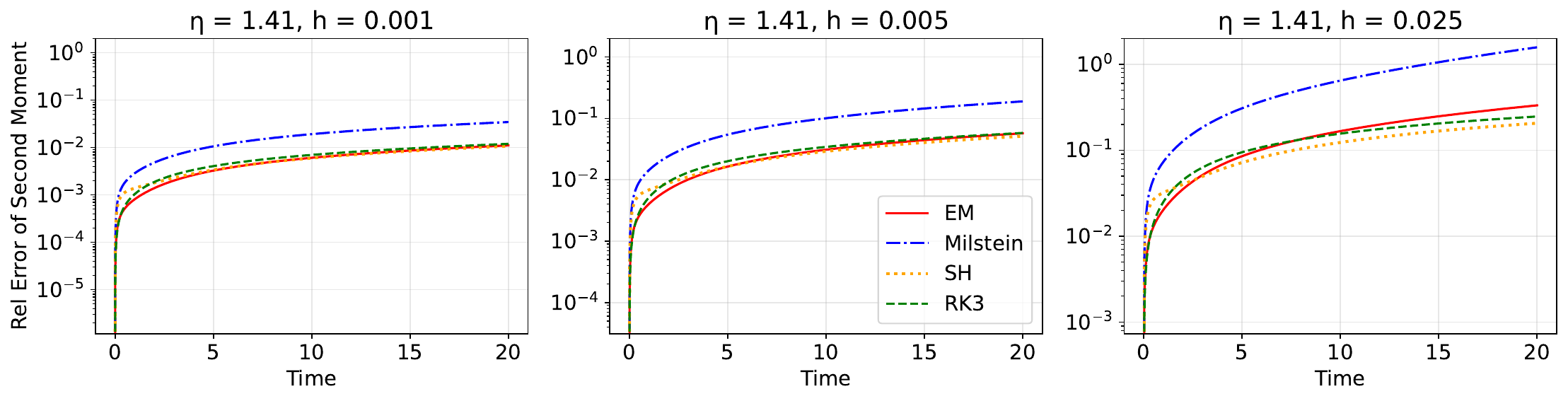}
    \caption{A comparison of the evolution of the relative error in second moments by the discretisation schemes versus the exact evolution up to time $T=20$. The first row depicts the case when $\eta=0.1$ for step sizes $h=0.001,0.005,0.025$, while the second row depicts the case when $\eta=1.41$ with the same range of step sizes. On long time-scales, we see substantial accuracy gains from higher-order schemes only when $\eta$ is small.}
    \label{fig:2ndMom}
\end{figure}

To investigate the evolution of error in the regime as the final time $T$ increases, in Figures~\ref{fig:1stMom} and \ref{fig:2ndMom}, we compare the analytic first and second moments computed using \eqref{eq:analytic1stmomentevolution}  and \eqref{eq:analytic2ndmomentevolution} with the corresponding discrete moments over time up to $T=20$. We note that due to our nondimensionalisation leading up to \eqref{eq:reducedSDE}, $T=1$ corresponds to a typical relaxation time-scale in the evolution of the equation. 

The numerical results agree with our rigorous analysis: in Figure~\ref{fig:1stMom}, we see that for $\eta=0.1$ and $\eta=1.41$, the EM and Mil schemes are asymptotically stable and accurate for the first moment, as their absolute errors decay exponentially as $T$ increases for the values of $\eta$ chosen and a range of time step sizes $h$. Meanwhile, the RK3 and SH schemes exhibit asymptotic bias, with their discrete first moments stabilising at a value which deviates from the analytical values on long time-scales. In Figure~\ref{fig:2ndMom}, we see similar agreement: second moments tend to approach the correct values with errors behaving linearly in the time step size $h$.

We note that focusing on time $T=1$, the higher-order schemes tend to perform better or at least comparably with the low-order schemes; this is to be expected from the existing weak-order results on fixed time intervals. However, when $T\geq 10$, we see significant gains in accuracy from using a higher-order scheme only for second moments and when $\eta$ is small. In agreement with our analytic results, these experiments therefore allow us to see that improved accuracy over fixed time intervals does not necessary translate to better accuracy in asymptotic statistics, and so the results in Table~\ref{tab:accuracy_summary} do not contradict the weak order results summarised in Table~\ref{tab:str_weak_accuracy_summary}.

\subsection{Numerical experiments in a nonlinear case}
\label{sec:nonlinear}
As a further illustration of the utility of the benchmark above, we also investigated a numerical example for a one-dimensional SDE with coefficients which are smooth and globally Lipschitz, but not affine, in the general form \eqref{eq:Ito_SDE}, testing whether the predictions made using the benchmark problem provide insight in this case.
One application area in which advection-diffusion equations for particle transport of this form arise naturally as a homogenised limit are in the study of flows through porous media where the statistical properties medium varies in space \cite{ito2012diffusion,ikeda2014stochastic}. In particular, variation in the size of pores through which fluid can move leads to coupled spatial variation in the flow velocity and molecular diffusion of the particles in the flow. At the single particle scale, this yields variation in the drift and diffusion coefficients, $f(x)$ and $g(x)$. Study of long-time statistics in this case is important for accumulation of pollutants in soils where groundwater is being pumped out of the water table.

The particular case we consider is a stochastic diffusion process on an infinite domain, governed by the SDE
\begin{equation}
\rd x_t = -A\,\mathrm{erf}\left(\frac{x_t-x^*}{B}\right)\rd t+ \left(C+D\,\mathrm{erf}\left(\frac{x_t-x^*}{E}\right)\right)\,\rd W_t,
\label{eq:nonlineartest}
\end{equation}
where we have chosen a drift $f(x) = -A\,\mathrm{erf}\big(\frac{x-x^*}{B}\big)$ and diffusion $g(x)=C+D\,\mathrm{erf}(\frac{x-x^*}{E})$. The coefficient functions exhibit sublinear growth and are globally Lipschitz, so strong solutions to the SDE exist.

Under the conditions that $A,B>0$ and $D>C>0$, it can be shown that this equation is ergodic with respect to an equilibrium distribution with exponentially decaying tails. Moreover, we can Taylor expand the coefficient functions about $x=x^*$ to obtain
\begin{align*}
    f(x) 
    &= f(x^*)+f'(x^*)(x-x^*)+O\left(|x-x^*|^2\right) 
    = \frac{-2 A}{B \sqrt{\pi}} (x-x^*)+O\left(|x-x^*|^2\right)\qquad\text{and}\\
    g(x)
    &=g(x^*)+g'(x^*)(x-x^*)+O\left(|x-x^*|^2\right)
    = C+\frac{2D}{E\sqrt{\pi}}(x-x^*) + O\left(|x-x^*|^2\right)
\end{align*}
Neglecting terms at quadratic and higher orders, we see that the corresponding linearised benchmark has
$$
\eta = \frac{g'(x^*)}{\sqrt{|f'(x^*)|}} = \sqrt{\frac{2BD^2}{AE^2\sqrt{\pi}}}.
$$

\begin{figure}[hp!]
\includegraphics[width = 0.94 \textwidth]{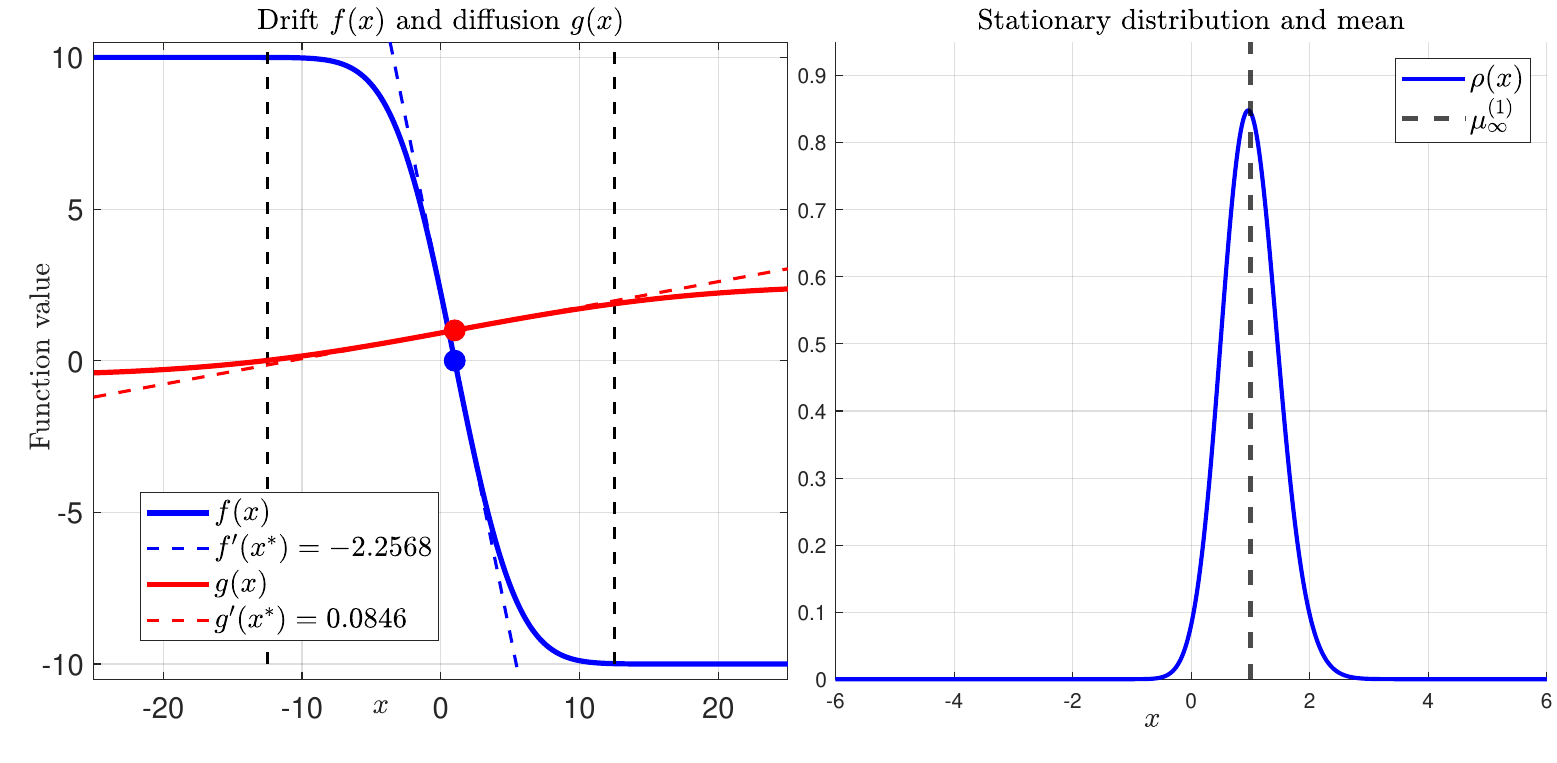}
\includegraphics[width = 0.94\textwidth]{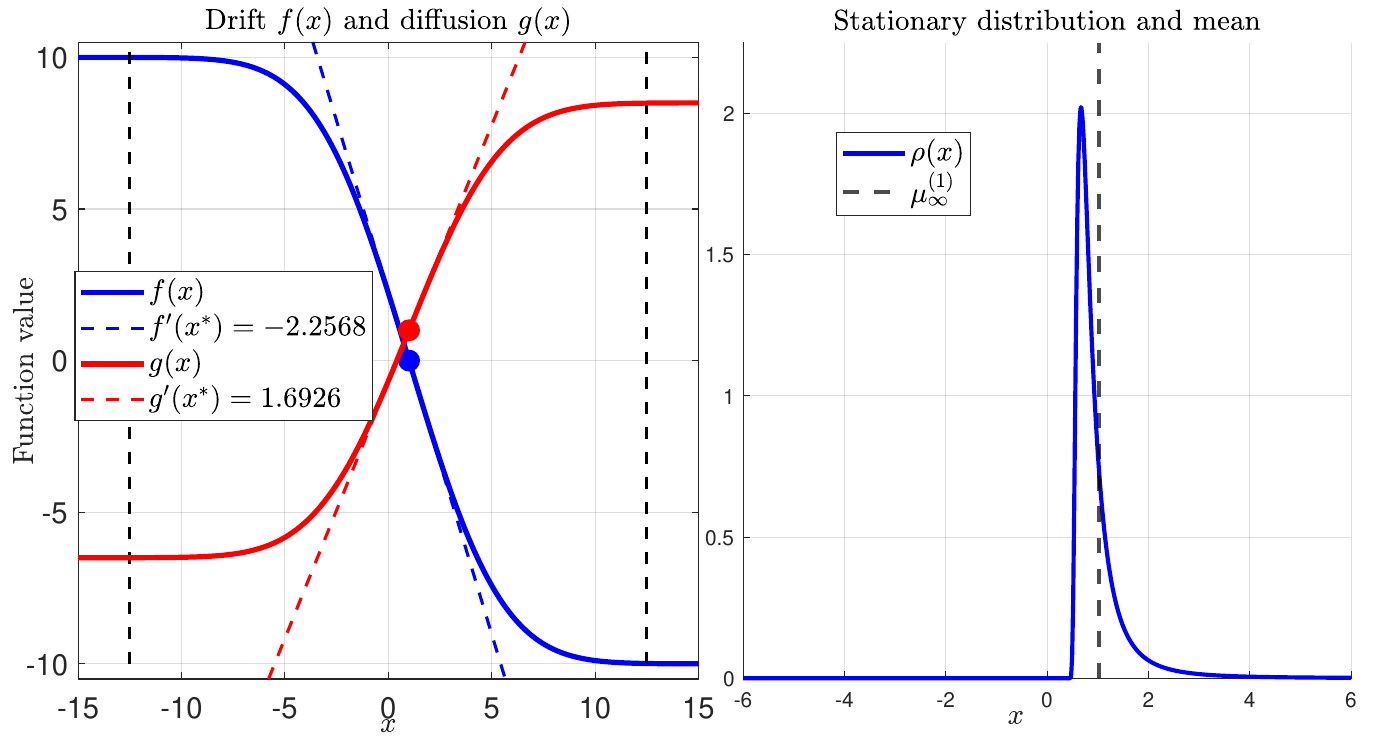}
\caption{Examples of the particular drift and diffusion coefficients used in the nonlinear example \eqref{eq:nonlineartest} (left column), along with the equilibrium measures computed by numerical integration (right column). The top row shows the case where $\eta\approx 0.05633$; while the bottom row shows the case where $\eta\approx 1.1267$. We see from the equilibrium measures that larger values of $\eta$ leads to significant skew in the distribution.}\label{Fig:porus_coef}
\end{figure}

For the purposes of numerical testing, we considered two sets of parameters, yielding different values of $\eta$:
\begin{enumerate}\itemsep2mm
\item $A =10, B=5, C= 1, D=1.5, E=20$, and $x^* = 1$, where $\eta\approx 0.05633$; and 
\item $A =10, B=5, C= 1, D=7.5, E=5$, and $x^* = 1$, where $\eta\approx  1.1267 $.
\end{enumerate}
The particular functions $f$ and $g$ in these cases are plotted in Figure~\ref{Fig:porus_coef}, along with the corresponding equilibrium distributions, which are computed numerically using the formula
\begin{equation}\label{nonlinear_exact_pdf}
p_\infty(x) = \frac{1}{Z g(x)^2}\exp\left(\int_{-\infty}^x\frac{2 f(y)}{g(y)^2}\rd y\right)
\qquad\text{where}\qquad Z:=\int_{-\infty}^\infty\frac{1}{g(x)^2}\exp\left(\int_{-\infty}^x\frac{2 f(y)}{g(y)^2}\rd y\right)\rd x.
\end{equation}
The density was computed using versions of the trapezoid rule over a uniform grid, using the functionality provided by the \texttt{scipy.integrate} package \cite{2020SciPy-NMeth}.

\subsubsection{Investigation of asymptotic moment accuracy}
\begin{figure}[htp!]
\centering 
\subfigure [$\eta\approx 0.05633$]
{\includegraphics[width = 0.46\linewidth, height = 6.75 cm]
{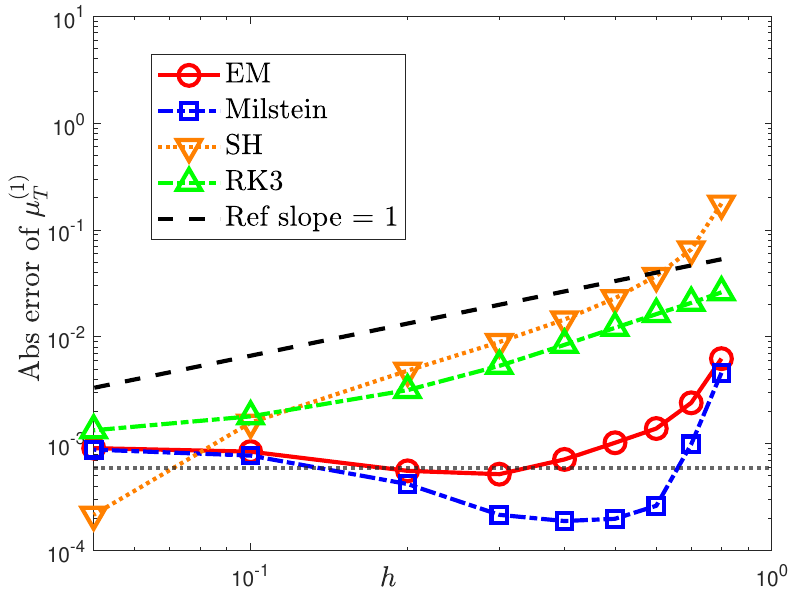}}
\subfigure [$\eta\approx 1.1267$]
{\includegraphics[width = 0.46\linewidth, height = 6.75 cm]
{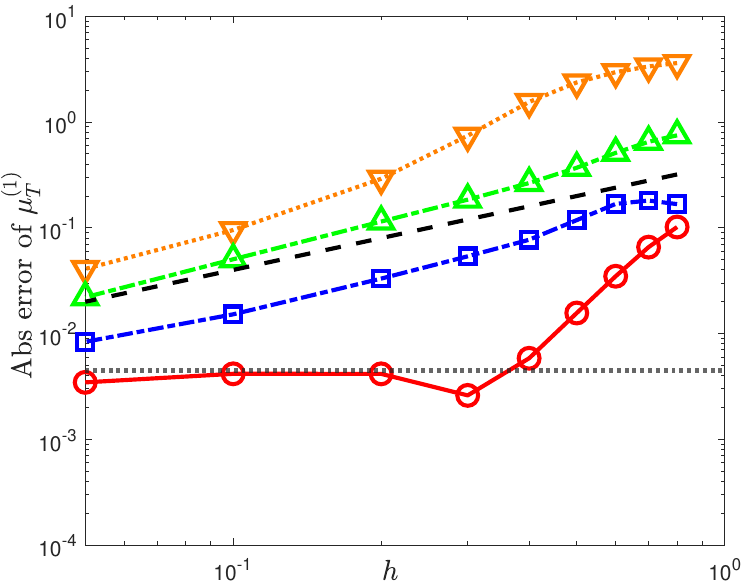}}
\caption{
Plots of the absolute error between the sample mean at time $T=20$ using different time-step sizes and schemes, and the stationary mean computed by numerical integration. Subfigure (a) shows the error in simulating \eqref{eq:nonlineartest} with parameter set (1), while (b) shows the results for parameter set (2). $6\times10^5$ sample trajectories were used for each point. The dotted line provides a reference value for the sampling error, computed as a 95\% confidence interval for the sample estimator with the number of samples used. In both test cases, the first-order schemes perform better than the higher-order schemes, with EM showing superlinear convergence up to the point where further reductions are limited by sampling error.} 
\label{Fig:1d_porousM_mean}
\end{figure}

To investigate whether our long-time asymptotic accuracy and stability results for the benchmark continue to be relevant for this nonlinear example, we simulated $N=6\times 10^5$ realizations of the dynamics of \eqref{eq:nonlineartest} using the various schemes analysed, and computed sample estimates of the first moment at time $T=20$. For both test cases, a range of step sizes were employed, with $h\in 0.05\times [1, 2:2:16]$. An analytic value for the equilibrium mean was not available for the processes considered, so a reference was computed by numerical integration. In this case, we do not have an exact value of the mean equilibrium distribution to use for comparison. Employing the trapezoid rule, the equilibrium distribution \eqref{nonlinear_exact_pdf} and its mean were computed and normalised on the interval $[-100,100]$ using a uniform grid with step size $\Delta x =0.0025$. Progressive refinement of the grid allowed showed this was sufficient to guarantee 5 digit accuracy in the reference values of asymptotic mean for the two test cases. The absolute differences between the sample mean and this reference are reported for each scheme and time-step size in Figure~\ref{Fig:1d_porousM_mean}.

To interpret these results, we note that the error between the simulated mean and the exact stationary mean arises from three sources:
\begin{enumerate}
    \item the finite-time truncation error, arising from taking the mean position at $T=20$, rather than $T=\infty$;
    \item the sampling error, arising from the use of a sample estimate of the mean with $N$ samples; and
    \item the discretization error which our results analyse, arising from the applying numerical scheme with $h>0$.
\end{enumerate}
For a final time $T=20$, the error in the mean arising from finite-time truncation with both sets of parameters is several orders of magnitude smaller than other sources, as the dynamics converges to equilibrium exponentially fast. This convergence can be observed directly in plots deriving from our later results (see Figures~\ref{Fig:fine_vs_coarse_small} and \ref{Fig:fine_vs_coarse_large}).
The sampling error can be assessed by applying the central limit theorem, allowing us to obtain a 95\% empirical confidence interval 
$
\pm z_{0.975} \frac{\sigma_N}{\sqrt{N_{\text{traj}}} }
$
about the sample mean, where $\sigma_N$ is a sample estimate of the standard deviation. For a fixed number of samples, the sampling error will limit the smallest possible error we can observe by sampling.

Plots of the absolute error between the sample mean and the reference value for the equilibrium mean obtained by numerical integration are plotted in Figure~\ref{Fig:1d_porousM_mean}. A reference value for the sampling error is shown as a dotted line, and we see that in line with the discussion above, for the schemes that reach this limit within the range tested, the error stabilises at around this value as the time-step shrinks. We see that the first-order scheme perform better in both test cases, with evidence of super-linear convergence for the EM scheme in both cases, and superlinear convergence for Mil in case (1).

Our observations do correlate strongly with the asymptotic moment accuracy analysis on the benchmark. Recall that it was shown that the asymptotic first moments are exact for EM and SH, and $O(h)$-accurate for SH and RK3. While we might not expect such results to carry over directly to this nonlinear setting, we do observe that the error in the mean predicted by the EM schemes behaves super-linearly for both test cases, and that both EM and Mil perform better on this quantity than the higher-order schemes, with sampling error saturating at much larger time-steps.

Our experiments here focus on first moments. Part of the reason for doing so is that the sample error for the statistical estimator of the second moment is controlled by the fourth moment of equilibrium distribution, and in both test cases, the fourth moment is large in these examples. This means that eliminating this source of sampling error to make a meaningful comparison among various schemes would have required excessive computational work, and was therefore not pursued.

\subsubsection{Investigation of moment stability}
Next, we investigate whether the predictions concerning the moment stability of the benchmark problem also contribute understanding of the stability properties in the nonlinear problems we consider.

First, we note that the stability regions shown in Figure~\ref{Fig:stab} are derived from the dimensionless form of the benchmark equation \eqref{eq:reducedSDE}. As such, to estimate the stability conditions for the nonlinear SDE \eqref{eq:nonlineartest}, which retains physical spatial and temporal units, a time rescaling is required. Specifically, the dimensionless time step $h_t$ in Figure~\ref{Fig:stab} must be mapped to the dimensional time step $h_{\tau}$ via the rescaling
\[
h_{\tau}\approx |f'(x^*)|^{-1} h_{t}. 
\]
Recall that it was shown that the EM and Mil schemes are asymptotically stable for the first moment for the dimensionless benchmark problem when $h_t<2$, so choosing a step-size corresponding to $h_t=2$ should reflect a point on the stability boundary for these schemes. Since the parameters controlling the deterministic drift $f$, $A$ and $B$, are fixed across the two test cases, we performed simulations for both test cases with an identical coarse dimensional time-step $H_\tau\approx 2|f'(x^*)|^{-1}\approx 0.8862$, and compared with a fine dimensional time-step $h_\tau=0.05$. Using the indicative values of $\eta$ given by the two case, we find that the SH scheme should be unstable at the coarser step-sizes in the two cases, and stable at the finer step-size in both cases. The RK3 scheme should be stable for both time-step sizes in case (1), and unstable for the larger step and stable for the finer step in case (2).

\begin{figure}[htp!]
\subfigure [$h=0.05, \quad \eta \approx 0.05633$]
{\includegraphics[width = 0.46\textwidth]{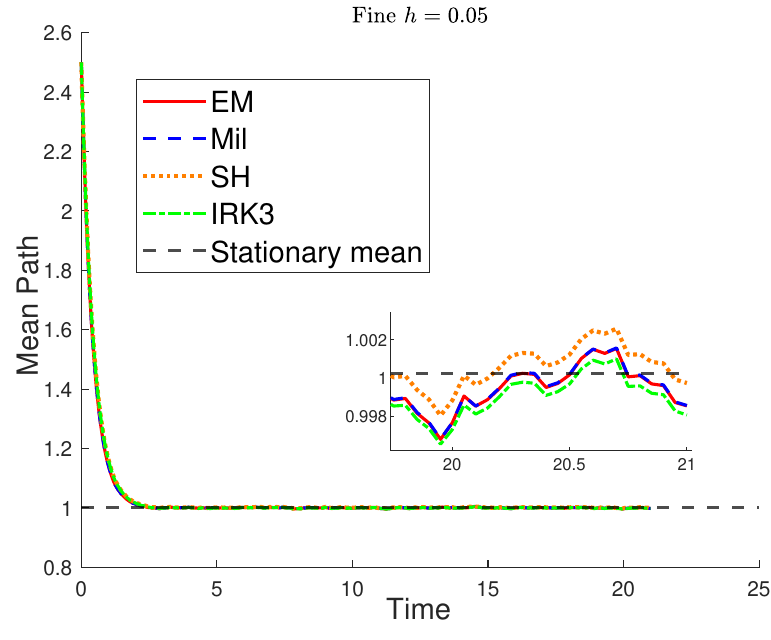}}
\subfigure [$H=0.8862,\quad \eta \approx 0.05633$]
{\includegraphics[width = 0.46\textwidth ]{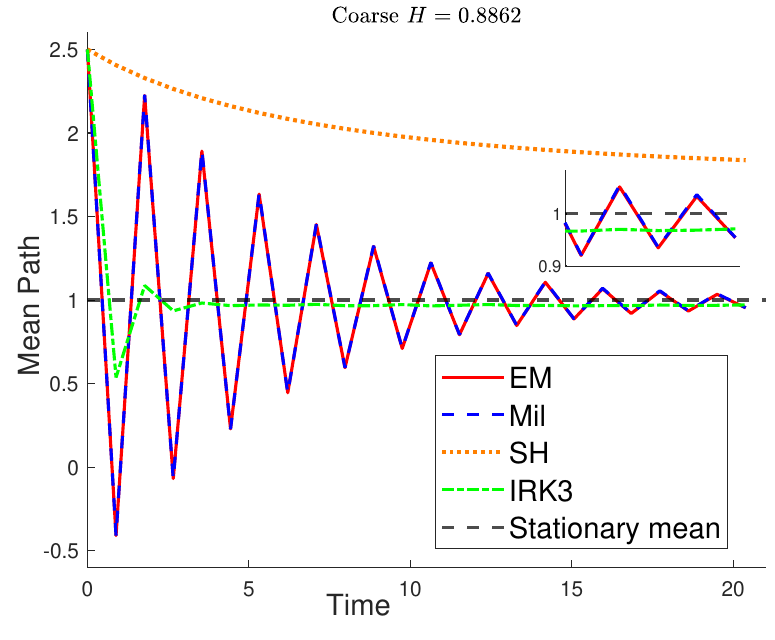}}
\caption{A comparison of the evolution of the first moment for discretisations of the nonlinear problem \eqref{eq:nonlineartest} in case (1), where parameters chosen predict a small value of $\eta\approx 0.05633$. The plots compare the means predicted with a fine step size $h = 0.05$ (left) and a coarse step size $H = 0.8882$ (right), which are averaged over $6\times10^5$ sample trajectories.
}\label{Fig:fine_vs_coarse_small}
\end{figure}

The results of the simulations in test case (1), where $\eta\approx 0.05633$ are shown in Figure~\ref{Fig:fine_vs_coarse_small}.
For fine step size $h=0.05$, all schemes are stable with mean paths approaching the stationary mean as predicted by the analysis. For the coarse step size $H=0.8862$, the sample means obtained from trajectories of EM and Mil schemes exhibit highly oscillatory behaviour, suggesting we are indeed close to the stability threshold. As predicted, the RK3 scheme is still stable with its numerical long-time mean approaching the analytical value. This is consistent with our benchmark analysis, as RK3 is the most stable scheme for small $\eta$ (see Fig~\ref{Fig:stab}). The SH scheme appears to equilibriate at a large biased prediction, which likely reflects instability of the scheme as predicted by the benchmark: the nature of the nonlinear system is such that no exponential growth can be expected (the forcing and noise tend to constants at large values).

\begin{figure}[htp!]
\subfigure [$h=0.05,\quad \eta\approx 1.1267$]
{\includegraphics[width = 0.46\textwidth]{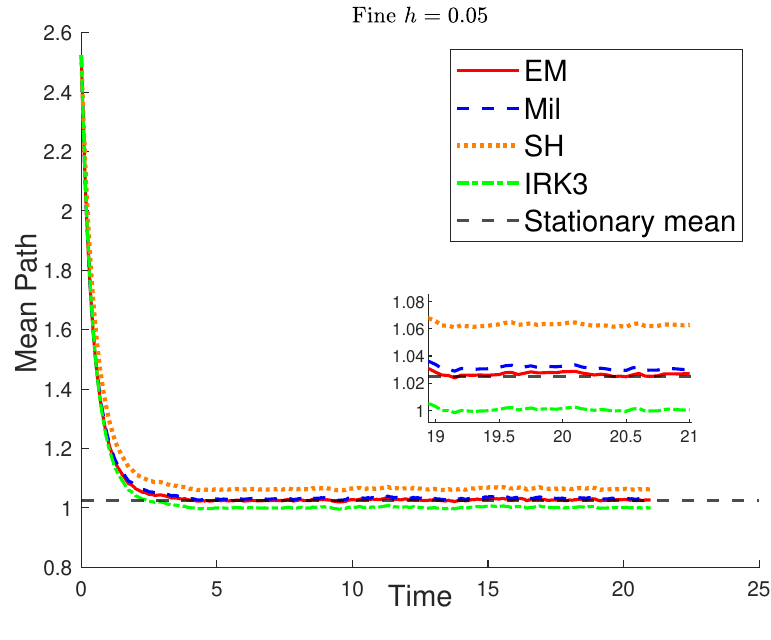}}
\subfigure [$H=0.8862,\quad \eta\approx 1.1267$]
{\includegraphics[width = 0.46\textwidth]{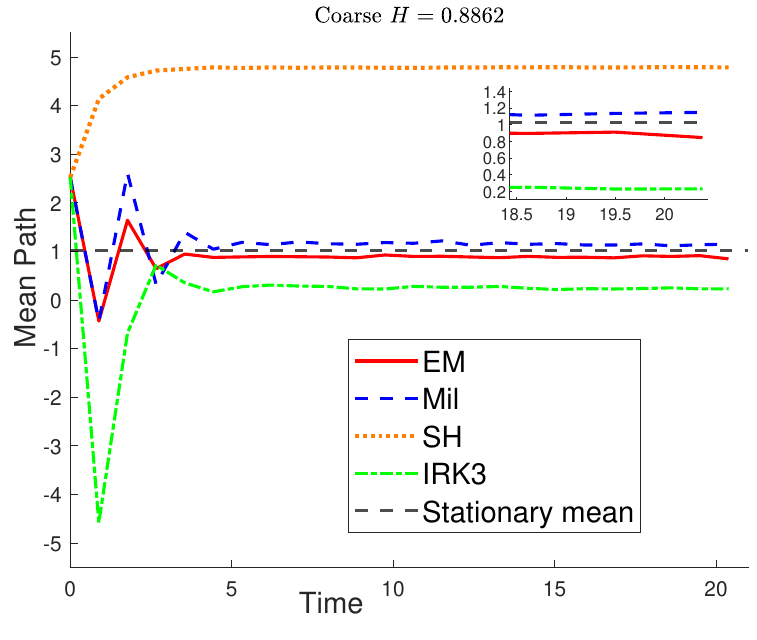}}
\caption{A comparison of the evolution of the first moment for discretisations of the nonlinear problem \eqref{eq:nonlineartest} in case (2), where the parameters chosen predict a large value of $\eta\approx 1.1267$. The plots compare the means predicted with a fine step size $h = 0.05$ (left) and a coarse step size $H = 0.8862$ (right), averaged over $6\times10^5$ sample trajectories.
}\label{Fig:fine_vs_coarse_large}
\end{figure}

The results of case (2) where $\eta\approx 1.1267$ are shown in Figure~\ref{Fig:fine_vs_coarse_large}. For the finer time-step size, we find that the lower-order schemes, EM and Mil, are more accurate in computing the long-time mean, as predicted by the analysis and as reflected by the accuracy plots in Figure~\ref{Fig:1d_porousM_mean}, with discrepancies visible for the higher-order schemes. For the coarse step-size $H=0.8862$, EM and Mil initially exhibit oscillations but later stabilize and produce long-time means with discrepancy less than $0.2$, even at this coarse step-size. In contrast, the RK3 exhibits dramatic initial oscillations and yields a long-time mean with a $0.8$ discrepancy from the analytical value. The SH scheme exhibits a growing error in the mean over time and appears to be unstable, consistent with the analysis. Again, there is no dramatic blow-up due to the nature of the system, but our stability analysis does indeed provide meaningful information regarding the performance of the schemes, even in the non-linear case.

The results for both values of $\eta$ with $h=0.05$ and $H=0.8862$ are consistent with the analysis, demonstrating that our study of the benchmark equation \eqref{eq:reducedSDE} remains relevant even when the coefficients are nonlinear.

\section{Conclusion}
Motivated by the need of practitioners to capture meaningful statistical averages of solutions to SDEs on long time-scales, we used a nondimensionalisation to derive a family of one-dimensional benchmark problems to evaluate the stability of numerical schemes for SDEs. We argued that in a generic case, the performance of numerical schemes on this benchmark provides more information about the long-time performance of these schemes for a range of applications than the more usual benchmark of geometric Brownian motion. We used our benchmark to present an analysis of the ability of four explicit numerical schemes used for simulating SDEs to recover first and second moments on long time-scales both accurately and stably. We observed that lower-order schemes such as the Euler-Maruyama tend to preserve asymptotic statistical accuracy of both first and second moments better than higher-order schemes, a surprising trend which persists even when testing a realistic nonlinear benchmark problem. Natural next steps would be to expand our analysis to include a range of implicit schemes, and to explore similar benchmarks for higher-order and higher-dimensional SDEs.

\subsection*{Acknowledgements}
The authors would like to thank Dr Josephine Evans at University of Warwick for her helpful discussions on the ergodic properties of SDEs during this work. We would also like to thank the anonymous referee for their detailed feedback which was of great help in improving the manuscript.

\subsection*{Data availability statement}
All the data and code required to reproduce the figures can be found in the Zenodo repository accompanying this work \cite{code}.

\subsection*{Open access statement}
For the purpose of open access, the authors have applied a Creative Commons Attribution (CC-BY) licence to any Author Accepted Manuscript
version arising from this submission.

\appendix
{
\section{Asymptotic moment stability for Geometric Brownian Motion}
\label{app:GBM}

In this appendix, we provide a parallel analysis of the stability of the schemes we study for geometric Brownian motion (GBM) as a benchmark.

\subsection{Evolution of moments for geometric Brownian motion}

\noindent
To provide a clear comparison between the previous GBM benchmark 
\[
\rd x_t = -x_t \rd_t + \eta x_t \rd W_t 
\]
used to define MS-stability
and our benchmark \eqref{eq:reducedSDE}, we compute the analytic and numerical first and second moments of various schemes for the GBM in this section. Taking the GBM from a deterministic initial position $x(0)=x_0$ and assuming that $\eta^2<2$, we have 
\begin{equation*}
         \rd x_t = -x_t\rd t+\eta x_t\rd W_t.
\end{equation*}
In this case, the Fokker-Planck equation for the probability density $p(x,t)$ of the state is given by 
\begin{equation*}
    \partial_tp=\partial_x(xp+\partial_x(\tfrac12\eta^2x^2p)).
\end{equation*}
As in Section~\ref{sec:moments}, we derive an equation for the first moment by taking time derivatives and integrating by parts, giving
\begin{equation*}
    \frac{d}{dt}\mu^{(1)}(t)=\int_{-\infty}^{\infty}x\partial_t p\rd x=\int_{-\infty}^{\infty}x\partial_x(xp)\,\rd x+\frac{\eta^2}{2}\int_{-\infty}^{\infty}x\partial_x^2(x^2p)\,\rd x=-\mu^{(1)}(t);
\end{equation*}
this equation has analytic solution $\mu^{(1)}(t)=\mu^{(1)}(0)e^{-t}$.
Similarly, the second moment satisfies
\begin{equation*}
    \frac{d}{dt}\mu^{(2)}(t)=\int_{-\infty}^{\infty}x^2\partial_tp\,\rd x = \int_{-\infty}^{\infty}x^2\partial_x(xp)\,\rd x+\frac{\eta^2}{2}\int_{-\infty}^{\infty}x^2\partial_x^2(x^2p)\rd x=-(2-\eta^2)\mu^{(2)}(t)
\end{equation*}
with solution $\mu^{(2)}(t)=\mu^{(2)}(0)e^{(-2+\eta^2)t}$. Since the initial condition is assumed to be deterministic, $p(x,0)=\delta_{x_0}(x)$, and we have the following analytical expressions and asymptotic values of the first and second moments:
\begin{equation*}
    \mu^{(1)}(t)=x_0\mathrm{e}^{-t}\to 0\quad\text{and}\quad
    \mu^{(2)}(t)= x_0^2\mathrm{e}^{-(2-\eta^2)t}\to 0\quad \text{as}\quad t\to \infty.    
\end{equation*}

\subsection{Asymptotic moment stability regions for various schemes}

In this section, we compute the asymptotic first and second-moment stability regions for the geometric Brownian motion under the four numerical schemes of interest.

\subsubsection{Euler Maruyama (EM)}
Taking expectations and using the scheme definition ~\eqref{eq:EM}, the first moment evolves according to
\begin{equation}\label{EM_1st_GBM}
    \mu_{n+1}^{(1)}=\mathbb{E}[x_{n+1}]= \mathbb{E}[(1-h)x_n+\eta x_nh^{1/2}\Delta W_n]=(1-h)\mu^{(1)}_n,
\end{equation}
and thus EM is first-moment asymptotically stable for $\vert 1-h\vert<1$, i.e., $0<h<2.$
Following a similar approach for the second moment, we have 
\begin{equation}\label{EM_2nd_GBM}
\begin{aligned}
    \mu_{n+1}^{(2)}&= \mathbb{E}[(x_n^2)(1-h+\eta h^{1/2}\Delta W_t)^2]\\
    &= \mathbb{E}[X_n^2]\left(
     (1-h)^2
     + 2\eta h^{1/2}(1-h)\mathbb{E}[\Delta W_n]
     + \eta^2h\mathbb{E}[(\Delta W_n)^2]
     \right) \\
     &= \mu_{n}^{(2)}\left(
     (1-h)^2
     + \eta^2 h
     \right),
\end{aligned}
\end{equation}
which leads to the stability requirement that 
\begin{equation*}
    0<h<2-\eta^2.
\end{equation*}
These stability requirements are identical to those found for the benchmark \eqref{eq:reducedSDE} in Section \ref{sec:stab}. In addition, the EM scheme for the GBM satisfies
\[
\text{EM:}\qquad 
\mu_{\infty}^{(1)}=0,\quad \mu_{\infty}^{(2)} =0 
\]
which show that EM scheme recovers the exact analytic values of asymptotic first and second moments for GBM.

\subsubsection{Milstein method (Mil)}

We see a similar phenomenon for the Mil method applied to the GBM. The Mil updates step for GBM is
\begin{equation*}
    X_{n+1} = X_n\left(1-h+\eta\Delta W_n+\frac{1}{2}\eta^2((\Delta W_n)^2-h)\right).
\end{equation*}
Taking the first moment, we have
\begin{equation}\label{Mil_1st_GBM}
    \mu_{n+1}^{(1)}= \mu_{n}^{(1)}(1-h).
\end{equation}
This leads to the stability condition for the first moment $0<h<2$.

\noindent
Taking the second moment, the recursion is
\begin{equation}\label{Mil_2nd_GBM}
    \mu_{n+1}^{(2)}= \mu_{n}^{(2)}\left(1+(\eta^2-2)h+(1+\frac{1}{2}\eta^4)h^2\right),
\end{equation}
which results in the asymptotic stability condition for the second moment to be
\begin{equation*}
    0<h<\frac{2-\eta^2}{1+\frac{1}{2}\eta^4}.
\end{equation*}
Again, the first and second stability conditions for GBM under Mil are identical to the stability conditions we found for the benchmark \eqref{eq:reducedSDE}. In addition for Mil, we have
\[
\text{Mil:}\qquad 
\mu_{\infty}^{(1)}=0,\quad \mu_{\infty}^{(2)} =0, 
\]
again, the Mil scheme recovers the exact asymptotic first and second moments without bias for GBM.

\subsubsection{Stochastic Heun (SH)}
The evolution of the first moment under SH reads
\begin{equation}\label{sH_1st_GBM}
    \mu_{n+1}^{(1)} = \left(1-h+\tfrac{1}{2}\big(1+\tfrac{1}{2}\eta^2\big)^2 h^2\right)\mu_n^{(1)} = \left(1-h+\tfrac{1}{2}h^2+\tfrac{1}{2}\eta^2 h^2+\tfrac{1}{8}\eta^4 h^2\right)\mu_n^{(1)},
\end{equation}
which gives the stability region 
\[
0 < h <
\frac{2}{\big(1+\frac{1}{2}\eta^2\big)^2}=\frac{8}{\eta^{4} + 4\eta^{2} + 4}
\]
Similarly, we have the recurrence relation for the second moment: 
$$
\begin{aligned}
    \mu_{n+1}^{(2)}&=\left(1\vphantom{\frac{1}{4}} + (\eta^2-2)h + \left(2-\eta^2-\tfrac{1}{4}\eta^4\right)h^2\right. \\
     &\qquad\qquad\left.\vphantom{\frac{1}{4}}+ \left(\tfrac{1}{4}\eta^6+\tfrac{3}{4}\eta^4-1\right)h^3 + \left( \tfrac{1}{64}\eta^8+\tfrac{1}{8}\eta^6+\tfrac{3}{8}\eta^4+\tfrac{1}{2}\eta^2+\tfrac{1}{4}\right) h^4\right)\mu_{n}^{(2)},
\end{aligned}\label{SH_2nd_GBM}
$$
which yields the stability condition
\[
\left| 1 + (\eta^2-2)h + \left(2-\eta^2-\tfrac{1}{4}\eta^4\right)h^2
     + \left(\tfrac{1}{4}\eta^6+\tfrac{3}{4}\eta^4-1\right)h^3 + \left( \tfrac{1}{64}\eta^8+\tfrac{1}{8}\eta^6+\tfrac{3}{8}\eta^4+\tfrac{1}{2}\eta^2+\tfrac{1}{4}\right) h^4\right|<1.
\]
These stability conditions for the first and second moments under SH are remain the same when compared to the results for the benchmark \eqref{eq:reducedSDE}. Moreover, we have the asymptotic result that
\[
\text{SH:}\qquad 
\mu_{\infty}^{(1)}=0,\quad \mu_{\infty}^{(2)} =0;
\]
so the exact analytic values of asymptotic moments are recovered without bias for GBM under SH discretisation scheme.

\subsubsection{3-stage Runge-Kutta (RK3)}
For RK3,  we have for the first moment
\begin{equation}\label{RK3_1st_GBM}
\begin{aligned}
\mu^{(1)}_{n+1}
&=\left(1
  - h
  + \left(\tfrac{1}{2} - \tfrac{1}{8}\eta^{4}\right) h^{2}
  - \left(\tfrac{1}{48}\eta^{6} + \tfrac{1}{8}\eta^{4}
         + \tfrac{1}{4}\eta^{2} + \tfrac{1}{6}\right) h^{3}\right)\mu^{(1)}_{n}\\
&= \left(1
    -h
    +\tfrac{1}{8} h^2 (4 - \eta^4)
    -\tfrac{1}{48} h^3 (2+\eta^2)^3\right)\,\mu^{(1)}_{n}.
\end{aligned}
\end{equation}
which indicates the first-moment is asymptotically stable when 
\begin{equation*}
   \left|1    -h    +\tfrac{1}{8} h^2 (4 - \eta^4)
    -\tfrac{1}{48} h^3 (2+\eta^2)^3\right|<1.
\end{equation*}
Next for the second moment, we have 
\[
\mu^{(2)}_{n+1}
= R^{(2)}(h,\eta)\,\mu^{(2)}_{n},
\]
where $R^{(2)}$ is a polynomial in $h$ of the form
\[
R^{(2)}(h,\eta)
= a_0 + a_1 h + a_2 h^2 + a_3 h^3 + a_4 h^4 + a_5 h^5 + a_6 h^6,
\]
with coefficients
\begin{align*}
a_0 &= 1, \quad 
a_1 = \eta^{2} - 2,\quad 
a_2 = \tfrac{1}{4}\eta^{4} - 2\eta^{2} + 2,\\[0.3em]
a_3 &= -\tfrac{1}{8}\eta^{6} + \tfrac{3}{2}\eta^{2} - \tfrac{4}{3},\quad 
a_4 = \tfrac{9}{64}\eta^{8} + \tfrac{7}{24}\eta^{6}
      - \tfrac{3}{8}\eta^{4} - \tfrac{1}{2}\eta^{2} + \tfrac{7}{12},\\[0.3em]
a_5 &= \tfrac{1}{48}\eta^{10} + \tfrac{5}{32}\eta^{8}
      + \tfrac{5}{12}\eta^{6} + \tfrac{5}{12}\eta^{4} - \tfrac{1}{6},\\[0.3em]
a_6 &= \tfrac{1}{2304}\eta^{12} + \tfrac{1}{192}\eta^{10}
      + \tfrac{5}{192}\eta^{8} + \tfrac{5}{72}\eta^{6}
      + \tfrac{5}{48}\eta^{4} + \tfrac{1}{12}\eta^{2} + \tfrac{1}{36}.
\end{align*}
Again, stability for the second moment requires that $|R^{(2)}(h,\eta)|<1$.
Both stability conditions for the first and second moments for GBM under RK3 are still the same compared to the benchmark \eqref{eq:reducedSDE}. Meanwhile, we have
\[
\text{RK3:}\qquad 
\mu_{\infty}^{(1)}=0,\quad \mu_{\infty}^{(2)} =0,
\]
which shows that RK3 recovers the exact second moment without bias for GBM.

\subsection{Summary}
For each of the schemes studied, the above computations show that the stability requirements for the GBM benchmark take the same form as for our benchmark equation \eqref{eq:reducedSDE} in terms of the parameter $\eta$ and timestep $h$. A stark difference however is that all schemes are asymptotically exact in terms of the first and second moment for GBM, with no bias. In part, this highlights the particularly special feature of GBM, which is that the drift and diffusion both vanish at the same point, and hence when it exists, the equilibrium measure is concentrated. As argued in Section~\ref{subsec:rescale}, this behaviour is rather special when considering an SDE with affine coefficients, highlighting the need for an alternative benchmark.}

\bibliographystyle{plain}
\bibliography{references}
\end{document}